\documentclass[12pt]{article}
\usepackage{bbm}
\usepackage{amsfonts}
\usepackage{pifont}
\usepackage{hyperref}

\usepackage{mathrsfs}
\usepackage{indentfirst}
\usepackage{amsmath}
\usepackage{amssymb}
\usepackage[amsmath, thmmarks]{ntheorem}
\usepackage{cite}
\usepackage{graphicx}
\usepackage{tikz}

\newtheorem{definition}{\bf Definition}[section]
\newtheorem{lemma}{\bf Lemma}[section]
\newtheorem{theorem}{\bf Theorem}[section]

\newtheorem{example}{\bf Example}[section]

\newtheorem{algorithm}{\bf Algorithm}[section]

\def\QEDopen{{\setlength{\fboxsep}{0pt}\setlength{\fboxrule}{0.2pt}\fbox{\rule[0pt]{0pt}{1.3ex}\rule[0pt]{1.3ex}{0pt}}}} 
\def\QED{\QEDopen}
\def\proof{{\bf Proof.} }
\def\endproof{\hspace*{\fill}~\QED\par\endtrivlist\unskip}

\begin{document}
\setcounter{page}{1}

\title{{\textbf{The best extending cover-preserving geometric lattices of semimodular lattices}}\thanks {Supported by the National Natural Science
Foundation of China (nos.11901064 and 12071325)}}
\author{Peng He$^1$\footnote{\emph{E-mail address}: hepeng@cuit.edu.cn}, Xue-ping Wang$^2$\footnote{Corresponding author}\\
\emph{\small 1. College of Applied Mathematics, Chengdu University of Information Technology}\\ \emph{\small Chengdu 610225, Sichuan, People's Republic of China}\\
\emph{\small 2. School of Mathematical Sciences, Sichuan Normal University}\\ \emph{\small Chengdu 610066, Sichuan, People's Republic of China}\\\emph{\small Corresponding author's e-mail address: xpwang1@hotmail.com}}
\newcommand{\pp}[2]{\frac{\partial #1}{\partial #2}}
\date{}
\maketitle

\begin{quote}
{\bf Abstract}

In 2010, G\'{a}bor Cz\'{e}dli and E. Tam\'{a}s Schmidt mentioned that the best cover-preserving embedding of a given semimodular lattice
is not known yet [A cover-preserving embedding of semimodular lattices into geometric lattices, Advances in Mathematics 225 (2010) 2455-2463]. That is to say: What are the geometric lattices $G$ such that a given finite semimodular lattice $L$ has a cover-preserving embedding into $G$ with the smallest $|G|$? In this paper, we propose an algorithm to calculate all the
best extending cover-preserving geometric lattices $G$ of a given semimodular lattice $L$ and prove that the length and the number of atoms of every best extending cover-preserving geometric lattice $G$ equal the length of $L$ and the number of non-zero join-irreducible elements of $L$, respectively. Therefore, we comprehend the best cover-preserving embedding of a given semimodular lattice.


\end{quote}

\section{Introduction}\label{intro}
Let $L$ be a lattice. For
all $a, b\in L$, $a\parallel b$ denotes that $a\ngeq b$ and
$a\nleq b$, and $a\nparallel b$ denotes that $a\geq b$ or
$a\leq b$. $a\prec b$ means that $a< b$ and there is no element $c\in L$
such that $a<c<b$, and $a\preceq b$ represents that $a\prec b$ or $a=b$.
The set of non-zero join-irreducible elements and
the set of atoms of $L$ will be denoted by
$J(L)$ and $A(L)$, respectively. The length of $L$, that is,
$\mbox{sup}\{n: L\mbox{ has an }(n+1)\mbox{-element chain}\}$,
will be denoted by $\ell(L)$.
Let $A$ and $B$ be two sets. We define $A-B=\{x\in A: x\notin B\}$. We
assume that the readers are familiar with the basic
notions of lattices such as a partially ordered set (poset), a chain, a
lattice, a distributive lattice, a modular
lattice, a semimodular lattice etc..
Here, we just recall a necessary concept from the theory of
lattices (see, e.g., \cite{Crawley73,Stern99}).  We say a lattice $L$ is (upper)
semimodular if $a\prec b$ implies $a\vee c\preceq b\vee c$ for all $a, b, c \in L$. We
know from Crawely and Dilworth \cite[Theorem 3.7]{Crawley73}
(see also \cite[Theorem 1.7.1]{Stern99}) that for
a strongly atomic algebraic lattice $L$, a semimodularity is equivalent to Birkhoff's condition:
\begin{equation*}
\mbox{ for all } a, b \in L, \mbox{ if } a\wedge b \prec a \mbox{ and } a\wedge b \prec b,
\mbox{ then } a\prec a\vee b \mbox{ and } b\prec a\vee b.
\end{equation*}
It is well known that if $L$ is a semimodular lattice with $\ell(L)=m$ and $|J(L)|=n$, then $n\geq m$ (see \cite{Stern99}).

Classically semimodular lattices arise out of certain closure operators satisfying what is now
usually called the Steinitz-Mac Lane exchange property. A semimodularity is one of the most important
links between combinatorics and lattice theory (see, e.g., \cite{Stern99, Czedli10}), and the structure
of a semimodular lattice plays an important role in lattice theory (see, e.g., \cite{Gratzer14, He}). A
particular interest is deserved by geometric lattices, originally called matroids, which are semimodular
atomistic lattices of finite length.

The Dilworth Embedding Theorem states that each finite lattice $L$ can be embedded in a finite geometric
lattice (see \cite{Crawley73}). Further, P. Pudl\'{a}k, J. T\"{u}ma\cite{Pudlak80} proved
that each finite lattice $L$ can be embedded in a finite partition lattice (finite partition lattices are
geometric lattices). In 1986, G. Gr\"{a}tzer and E. W. Kiss\cite{Gratzer86} showed that each finite
semimodular lattice $L$ has a cover-preserving embedding into a finite geometric lattice. Recently, G. Cz\'{e}dli and E. T. Schmidt \cite{Czedli10} extended the results in \cite{Gratzer86}, and they proved that each semimodular lattice $L$ of finite length
has a cover-preserving embedding into a geometric lattice $G$ of the same length and the number
of atoms of $G$ equals the number of non-zero join-irreducible elements of $L$. That is, they proved the following theorem.
\begin{theorem}[\cite{Czedli10}]\label{th1}
Let $L$ be a semimodular lattice of finite length. Then there exists a geometric
lattice $G$ such that $L$ is a cover-preserving sublattice of $G$, $|J(L)|=|A(G)|$ and $\ell(L)=\ell(G)$.
\end{theorem}
Finally, they mentioned that the best cover-preserving embedding
is not known yet. That is to say: What are the geometric lattices $G$ such that a given finite semimodular lattice $L$ has a cover-preserving embedding into $G$ with the smallest $|G|$? In this paper, we shall construct all the best
cover-preserving embeddings of a given finite semimodular lattice $L$ into
geometric lattices $G$ and prove that the length and the number of atoms of every best extending cover-preserving geometric lattice $G$ equal the length of $L$ and the number of non-zero join-irreducible elements of $L$, respectively.

For the detailed information on semimodular lattices and partially ordered sets
the readers are referred to  \cite{Birkhoff73,Crawley73,Gratzer14,Stern99}.
We use the terminologies and notations of \cite{Crawley73, Birkhoff73}.

\section{Atomistic partially ordered sets}

In this section, we shall introduce the concept of an atomistic partially ordered set
and then investigate some of its basic properties.

\begin{definition}\label{d21}
\emph{Let $(P,\leq)$ be a finite partially ordered set and
$$\ell(P)=\mbox{sup}\{n: P\mbox{ has an }(n+1)\mbox{-element chain}\}. $$
Then we say that $\ell(P)$ is the length of $P$.}

\emph{If $P$ has the minimum element $0$, then let $\ell_P(x)$, or $\ell(x)$ for brevity,
denote the length of $[0,x]$ for each element $x\in P$. Thus, $\ell(0)=0$ and $\ell(1)=\ell(P)$ when $1$ is
the maximum element of $P$.}
\end{definition}

Similar to the definitions of atoms of lattices, an element that covers the least element $0$ of a
partially ordered set $P$ is referred to as an atom of $P$, and denoted by $A(P)$ the set of
atoms of $P$, i.e., $A(P)=\{x\in P: x\succ 0\}$. In particular, $A_P(y)=A([0,y])$, or $A(y)=A([0,y])$ for brevity, for each $y\in P$.

\begin{example}\label{e21}
\emph{The Hasse diagram of a partially ordered set $P$ is shown as Fig.1.}
 \end{example}
\par\noindent\vskip50pt
 \begin{minipage}{11pc}

\setlength{\unitlength}{0.75pt}\begin{picture}(600,200)

\put(170,40){\circle{6}}\put(168,28){\makebox(0,0)[l]
{\footnotesize $0$}}

\put(100,110){\circle{6}}\put(85,105){\makebox(0,0)[l]{\footnotesize$
  a$}}

  \put(380,110){\circle{6}}\put(385,105){\makebox(0,0)[l]{\footnotesize$
    c$}}

   \put(240,110){\circle{6}}\put(255,110){\makebox(0,0)[r]{\footnotesize$
   b$}}

   \put(100,180){\circle{6}}\put(85,180){\makebox(0,0)[l]{\footnotesize$
   x$}}

   \put(240,180){\circle{6}}\put(255,180){\makebox(0,0)[r]{\footnotesize$
      y$}}

      \put(170,250){\circle{6}}\put(175,260){\makebox(0,0)[r]{\footnotesize$
   1$}}

  \put(100,114){\line(0,1){62}}
  \put(168,43.5){\line(-1,1){64.5}}
  \put(173,42.5){\line(1,1){64.5}}
  \put(173,41.5){\line(3,1){203}}

 \put(377,112.5){\line(-3,2){203}}
  \put(240,114){\line(0,1){62}}

  \put(103,112.5){\line(2,1){133}}
  \put(102,183.5){\line(1,1){64.5}}
  \put(237,112.5){\line(-2,1){133}}
  \put(238,183.5){\line(-1,1){64.5}}
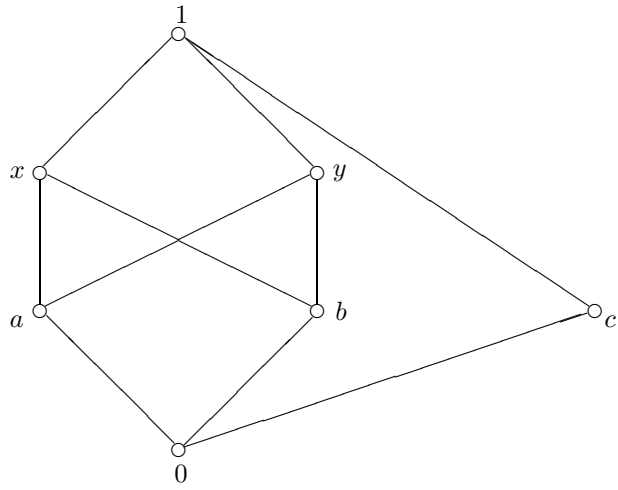
    \put(105,5){ Fig.1 The partially ordered set $P$.}
  \end{picture}
  \end{minipage}\\
In Fig.1, $A(P)=A(1)=\{a, b, c\}$, $A(x)=A(y)=\{a, b\}$ and $A(0)=\emptyset$.

\begin{definition}\label{d22}
\emph{A finite partially ordered set $P$ with the minimum element $0$ is atomistic if and only if
the following two conditions are satisfied: for all $x, y\in P$,\\
(1) $x< y$ implies that $A(x)\subsetneq A(y)$;\\
(2) $x\parallel y$ yields that $A(x)\nsubseteq A(y)$ and $A(y)\nsubseteq A(x)$.}
\end{definition}

\par\noindent\vskip50pt
 \begin{minipage}{11pc}

\setlength{\unitlength}{0.75pt}\begin{picture}(600,280)

\put(205,40){\circle{6}}\put(200,28){\makebox(0,0)[l]
{\footnotesize $0_P$}}

\put(100,110){\circle{6}}\put(85,110){\makebox(0,0)[l]{\footnotesize$
  1$}}

  \put(170,110){\circle{6}}\put(177,110){\makebox(0,0)[l]{\footnotesize$
   2$}}

   \put(240,110){\circle{6}}\put(252,110){\makebox(0,0)[r]{\footnotesize$
   3$}}

\put(310,110){\circle{6}}\put(322,110){\makebox(0,0)[r]{\footnotesize$
  4$}}

   \put(205,180){\circle{6}}

   \put(170,250){\circle{6}}

      \put(240,250){\circle{6}}

  \put(380,250){\circle{6}}

    \put(240,320){\circle{6}}\put(240,330){\makebox(0,0)[l]
{\footnotesize $1_P$}}

  \put(207,43.5){\line(1,2){31.3}}
  \put(203,43.5){\line(-1,2){31.3}}
  \put(207,43.5){\line(3,2){99.3}}
  \put(203,43.5){\line(-3,2){99.3}}
  \put(172,113.5){\line(1,2){31.3}}
  \put(238,113.5){\line(-1,2){31.3}}
  \put(207,183.5){\line(1,2){31.3}}
  \put(203,183.5){\line(-1,2){31.3}}
   \put(240,254){\line(0,1){62}}
    \put(102,113.5){\line(1,2){66.3}}
    \put(312,113.5){\line(1,2){66.3}}
    \put(308,113.5){\line(-1,2){66.3}}
     \put(172,253.5){\line(1,1){64.3}}
    \put(378,253.5){\line(-2,1){134.3}}
    \put(102,113.5){\line(2,1){273.3}}
    \put(172,113.5){\line(3,2){204.3}}
    \put(105,5){ Fig.2 The atomistic partially ordered set $P$.}
  \end{picture}
  \end{minipage}\\

By Definition \ref{d22}, one can check that Fig.1 is not atomistic since $x\parallel y$
but $A(x)=A(y)$, and Fig.2 is atomistic. Clearly, every finite atomistic lattice
is an atomistic partially ordered set, but the inverse is not true generally. For example,
Fig.2 is an atomistic partially ordered set, but it is not a finite atomistic lattice since it is not a lattice. However, the following lemma is clear.

\begin{lemma}\label{l21}
If a finite atomistic partially ordered set $P$ is a lattice, then $P$ is an atomistic lattice.
\end{lemma}

Let $\mathcal{P}(X)$ be the power set of a nonempty set $X$. Then we easily verify the following lemma.

\begin{lemma}\label{t21}
Let $|X|<\infty$ and $P\subseteq \mathcal{P}(X)$. If $\emptyset \in P$
and $\{\{x\}: x\in X\}\subseteq P$ then $(P, \subseteq)$ is a finite atomistic partially ordered set.
\end{lemma}

For convenience, in the following, if $P$ is a finite atomistic partially ordered set then we denote $\mathcal{S}_P=\{A(x): x\in P\}$.
\begin{lemma}\label{t22}
If $P$ is a finite atomistic partially ordered set, then
$(P, \leq)\cong(\mathcal{S}_P, \subseteq)$.
\end{lemma}
\proof
For $x\in P$, define $f: P\longrightarrow \mathcal{S}_P$ to be a map such that
$$f(x)=A(x) \mbox{ for any } x\in P.$$  We will show that the map $f$ is an isomorphism
of partially ordered sets.

It is clear that the map $f$ is well-defined. If $x, y\in P$ and $x\neq y$, then $f(x)=A(x)\neq A(y)=f(y)$ by Definition \ref{d22}. Hence, the map $f$ is injective. Moreover, it is clearly that
there exists $x\in P$ such that $U=A(x)=f(x)$ for any $U\in \mathcal{S}_P$ from the definition of $\mathcal{S}_P$, i.e., the map $f$ is surjective.
Consequently, the map $f$ is a one-to-one map. Below, we only need to prove that
both $f$ and its inverse are order-preserving.

Set $x, y\in P$ and $x< y$, and observe that application of the condition (1) of Definition \ref{d22} yields $f(x)=A(x)\subsetneq A(y)=f(y)$. Thus the map $f$ is order-preserving. Now suppose that $U, V\in \mathcal{S}_P$
and $U\subsetneq V$. Then there exist $x, y\in P$ such that $U=A(x)\subsetneq V=A(y)$. By
Definition \ref{d22}, $x< y$. Thus, the inverse of $f$ is order-preserving. Therefore,
$(P, \leq)\cong(\mathcal{S}_P, \subseteq)$.
\endproof

By Lemma \ref{t22}, every finite atomistic partially ordered set can be
considered as a set of sets. For instance, Fig.2 and Fig.3 are isomorphic.

\par\noindent\vskip50pt
 \begin{minipage}{11pc}

\setlength{\unitlength}{0.75pt}\begin{picture}(600,300)

\put(205,40){\circle{6}}\put(202,28){\makebox(0,0)[l]
{\footnotesize $\emptyset$}}

\put(100,110){\circle{6}}\put(75,105){\makebox(0,0)[l]{\footnotesize$
  \{1\}$}}

  \put(170,110){\circle{6}}\put(175,106){\makebox(0,0)[l]{\footnotesize$
    \{2\}$}}

   \put(240,110){\circle{6}}\put(265,110){\makebox(0,0)[r]{\footnotesize$
   \{3\}$}}

\put(310,110){\circle{6}}\put(333,110){\makebox(0,0)[r]{\footnotesize$
   \{4\}$}}

   \put(205,180){\circle{6}}\put(165,180){\makebox(0,0)[l]{\footnotesize$
   \{2,3\}$}}

   \put(170,250){\circle{6}}\put(160,250){\makebox(0,0)[r]{\footnotesize$
      \{1,2,3\}$}}

      \put(240,250){\circle{6}}\put(290,250){\makebox(0,0)[r]{\footnotesize$
   \{2,3,4\}$}}

  \put(380,250){\circle{6}}\put(430,250){\makebox(0,0)[r]{\footnotesize$
   \{1,2,4\}$}}

    \put(240,320){\circle{6}}\put(266,333){\makebox(0,0)[r]{\footnotesize$
   \{1,2,3,4\}$}}

  \put(207,43.5){\line(1,2){31.3}}
  \put(203,43.5){\line(-1,2){31.3}}
  \put(207,43.5){\line(3,2){99.3}}
  \put(203,43.5){\line(-3,2){99.3}}
  \put(172,113.5){\line(1,2){31.3}}
  \put(238,113.5){\line(-1,2){31.3}}
  \put(207,183.5){\line(1,2){31.3}}
  \put(203,183.5){\line(-1,2){31.3}}
   \put(240,254){\line(0,1){62}}
    \put(102,113.5){\line(1,2){66.3}}
    \put(312,113.5){\line(1,2){66.3}}
    \put(308,113.5){\line(-1,2){66.3}}
     \put(172,253.5){\line(1,1){64.3}}
    \put(378,253.5){\line(-2,1){134.3}}
    \put(102,113.5){\line(2,1){273.3}}
    \put(172,113.5){\line(3,2){204.3}}
    \put(100,5){ Fig.3 The atomistic partially ordered set $\mathcal{S}_P$.}
  \end{picture}
  \end{minipage}

\begin{definition}\label{d23}
\emph{Let $P$ be a finite atomistic partially ordered
set. A map $\mathfrak{I}_P$ from $P$ to the power sets of $\mathcal{P}(A(P))$ is called an independent function on $P$ if it has the following two properties: for any $x\in P$,\\
(1) if $\ell(x)=0$, then $\mathfrak{I}_P(x)=\{\emptyset\}$;\\
(2) if $\ell(x)\geq 1$, then
$$\mathfrak{I}_P(x)=\{S\cup\{a\}: S\in \mathfrak{I}_P(y), a\in A(x)-A(y), \ell(x)=\ell(y)+1 \mbox{ and }x\succ y\}.$$}
\end{definition}

Clearly, $\bigcup \mathfrak{I}_P(x)=A(x)$ for any $x\in P$. Let $P$ be a finite atomistic lattice. If $x, y\in P$, $x\parallel y$, then $\sigma \nsubseteq A(y)$ for any $\sigma\in \mathfrak{I}_P(x)$.

From Definition \ref{d23} and
Theorem 6.5 in \cite{Crawley73}, the following lemma is obviously.

\begin{lemma}\label{l22}
Let $L$ be a finite geometric lattice and $x\in L$. Then the following three
conditions are equivalent:\\
\emph{(1)} $S\in \mathfrak{I}_L(x)$;\\
\emph{(2)} $S$ is a maximal independent set of atoms of $[0, x]$;\\
\emph{(3)} $S$ is an independent set of atoms of $L$ and $\bigvee S=x$.
\end{lemma}

The diamond $M_3$ (see Fig.4) is a geometric lattice and $\mathfrak{I}_{M_3}(1)=\{\{a,b\},\{a,c\},\{b,c\}\}$.
One can verify that $a\vee b=a\vee c=b\vee c$=1, and $\{a,b\},\{a,c\}$ and $\{b,c\}$ are maximal
independent sets of atoms of $M_3$.

\par\noindent\vskip50pt
\begin{minipage}{11pc}
\setlength{\unitlength}{0.75pt}\begin{picture}(600,167)
\put(240,40){\circle{6}}\put(237,28){\makebox(0,0)[l]
{\footnotesize $0$}}
\put(240,200){\circle{6}}\put(238,209){\makebox(0,0)[l]{\footnotesize$
 1$}}
\put(160,120){\circle{6}}\put(145,120){\makebox(0,0)[l]
{\footnotesize $a$}}
\put(240,120){\circle{6}}\put(250,120){\makebox(0,0)[l]
{\footnotesize $b$}}
\put(320,120){\circle{6}}\put(330,120){\makebox(0,0)[l]
{\footnotesize $c$}}

\put(240,43){\line(0,1){74}}
\put(242,42){\line(1,1){76}}

\put(240,123){\line(0,1){74}}

\put(238,42){\line(-1,1){76}}

\put(162,122){\line(1,1){76}}

\put(318,122){\line(-1,1){76}}
\put(200,5){ Fig.4 $M_3$.}
\end{picture}
\end{minipage}

\section{Constructions of geometric lattices}
For the rest of this paper, unless otherwise stated, let $L$ be a fixed finite semimodular lattice. Following the convention
of, say, Crawley and Dilworth \cite{Crawley73} or Birkhoff \cite{Birkhoff73}, we assume that $L$ is
non-empty. Let $H(L)=L-A(L)\cup \{0\}$. For $x\in H(L)$, let $\Delta(x)$ be
a finite set satisfying that $\Delta(x)\cap L=\emptyset$ and $\Delta(x)\cap \Delta(y)=\emptyset$ while $x\neq y$, where $\Delta(x)$ may be empty set.
Insert every element in $\Delta(x)$ into $L$.
Extend the original order by $0\prec x^{'}\prec x$ for every $x^{'}\in \Delta(x)$;
this way we obtain a finite partially ordered set $(P, \leq)$ with $P=L\cup \bigcup_{x\in H(L)}\Delta(x)$. Notice that if $(P, \leq)$ is a lattice, then we call it an extending lattice of $L$.
The constructions of three new finite partially ordered sets
$P_1, P_2$ and $P_3$ are depicted in Fig.5; the black-filled elements are the inserted ones.

\par\noindent\vskip50pt
\begin{minipage}{11pc}
\setlength{\unitlength}{0.75pt}\begin{picture}(600,210)
\put(40,40){\circle{6}}\put(36,28){\makebox(0,0)[l]
{\footnotesize $0$}}
\put(40,80){\circle{6}}\put(25,78){\makebox(0,0)[l]
{\footnotesize$a$}}
\put(40,120){\circle{6}}\put(45,118){\makebox(0,0)[l]
{\footnotesize $b$}}
\put(40,200){\circle{6}}\put(38,210){\makebox(0,0)[l]
{\footnotesize $1$}}
\put(80,160){\circle{6}}\put(85,158){\makebox(0,0)[l]
{\footnotesize $c$}}
\put(0,160){\circle{6}}\put(-15,158){\makebox(0,0)[l]
{\footnotesize $d$}}

\put(40,43){\line(0,1){34}}
\put(40,83){\line(0,1){34}}
\put(42,122){\line(1,1){36}}
\put(38,122){\line(-1,1){36}}
\put(2,162){\line(1,1){36}}
\put(78,162){\line(-1,1){36}}
\put(30,5){$L$}

\put(170,40){\circle{6}}\put(168,28){\makebox(0,0)[l]
{\footnotesize $0$}}
\put(170,80){\circle{6}}\put(155,78){\makebox(0,0)[l]
{\footnotesize$a$}}
\put(170,120){\circle{6}}\put(175,118){\makebox(0,0)[l]
{\footnotesize $b$}}
\put(170,200){\circle{6}}\put(168,210){\makebox(0,0)[l]
{\footnotesize $1$}}
\put(210,160){\circle{6}}\put(215,158){\makebox(0,0)[l]
{\footnotesize $c$}}
\put(130,160){\circle{6}}\put(115,158){\makebox(0,0)[l]
{\footnotesize $d$}}
\put(130,80){\circle{6}}\put(115,78){\makebox(0,0)[l]
{\footnotesize $b^{'}$}}
\put(126,76){$\bullet$}
\put(170,43){\line(0,1){34}}
\put(170,83){\line(0,1){34}}
\put(172,122){\line(1,1){36}}
\put(168,122){\line(-1,1){36}}
\put(168,122){\line(-1,1){36}}
\put(132,162){\line(1,1){36}}
\put(208,162){\line(-1,1){36}}
\put(168,42){\line(-1,1){36}}
\put(132,82){\line(1,1){36}}
\put(166,5){$P_1$}

\put(340,40){\circle{6}}\put(336,28){\makebox(0,0)[l]
{\footnotesize $0$}}
\put(340,80){\circle{6}}\put(325,78){\makebox(0,0)[l]
{\footnotesize$a$}}
\put(340,120){\circle{6}}\put(325,118){\makebox(0,0)[l]
{\footnotesize $b$}}
\put(340,200){\circle{6}}\put(338,210){\makebox(0,0)[l]
{\footnotesize $1$}}
\put(380,160){\circle{6}}\put(385,158){\makebox(0,0)[l]
{\footnotesize $c$}}
\put(300,160){\circle{6}}\put(285,158){\makebox(0,0)[l]
{\footnotesize $d$}}
\put(296,76){$\bullet$}
\put(300,80){\circle{6}}\put(285,78){\makebox(0,0)[l]
{\footnotesize $b^{'}$}}
\put(296,76){$\bullet$}
\put(260,80){\circle{6}}\put(245,78){\makebox(0,0)[l]
{\footnotesize $d^{'}$}}
\put(256,76){$\bullet$}
\put(380,80){\circle{6}}\put(385,78){\makebox(0,0)[l]
{\footnotesize $c^{'}$}}
\put(376,76){$\bullet$}
\put(340,43){\line(0,1){34}}
\put(340,83){\line(0,1){34}}
\put(342,122){\line(1,1){36}}
\put(338,122){\line(-1,1){36}}
\put(338,122){\line(-1,1){36}}
\put(302,162){\line(1,1){36}}
\put(378,162){\line(-1,1){36}}
\put(338,42){\line(-1,1){36}}
\put(302,82){\line(1,1){36}}
\put(338,42){\line(-2,1){74.5}}
\put(262,82){\line(1,2){37.5}}
\put(342,42){\line(1,1){36}}
\put(380,83){\line(0,1){74}}
\put(333,5){$P_2$}

\put(500,40){\circle{6}}\put(499,28){\makebox(0,0)[l]
{\footnotesize $0$}}
\put(500,80){\circle{6}}\put(505,78){\makebox(0,0)[l]
{\footnotesize$a$}}
\put(500,120){\circle{6}}\put(505,118){\makebox(0,0)[l]
{\footnotesize $b$}}
\put(500,200){\circle{6}}\put(498,210){\makebox(0,0)[l]
{\footnotesize $1$}}
\put(540,160){\circle{6}}\put(545,158){\makebox(0,0)[l]
{\footnotesize $c$}}
\put(460,160){\circle{6}}\put(445,158){\makebox(0,0)[l]
{\footnotesize $d$}}
\put(456,76){$\bullet$}
\put(460,80){\circle{6}}\put(445,90){\makebox(0,0)[l]
{\footnotesize $d_1^{'}$}}
\put(416,76){$\bullet$}
\put(420,80){\circle{6}}\put(415,92){\makebox(0,0)[l]
{\footnotesize $d^{'}$}}
\put(536,76){$\bullet$}
\put(540,80){\circle{6}}\put(545,78){\makebox(0,0)[l]
{\footnotesize $c^{'}$}}
\put(476,76){$\bullet$}
\put(480,80){\circle{6}}\put(475,92){\makebox(0,0)[l]
{\footnotesize $b^{'}$}}
\put(500,43){\line(0,1){34}}
\put(500,83){\line(0,1){34}}
\put(502,122){\line(1,1){36}}
\put(498,122){\line(-1,1){36}}
\put(498,122){\line(-1,1){36}}
\put(462,162){\line(1,1){36}}
\put(538,162){\line(-1,1){36}}
\put(498,42){\line(-1,1){36}}
\put(460,83){\line(0,1){74}}
\put(498,42){\line(-2,1){74.5}}
\put(422,82){\line(1,2){37.5}}
\put(502,42){\line(1,1){36}}
\put(540,83){\line(0,1){74}}
\put(482,82){\line(1,2){17.5}}
\put(498,42){\line(-1,2){17.5}}
\put(496,5){$P_3$}
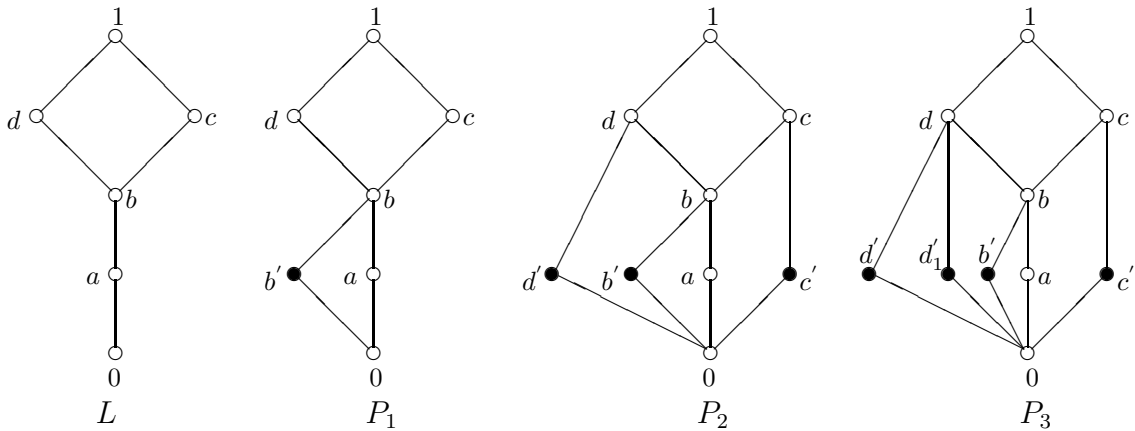
\put(5,-20){ Fig.5 An example of $L$ and the three extensions $P_1, P_2$ and $P_3$, respectively.}
\end{picture}
\end{minipage}

$ $\\

\begin{definition}\label{d31}
\emph{If $\Delta(x)\neq \emptyset$ for
every $x\in J(L)\cap H(L)$, then $(P, \leq)$
is called an extending standard form of $L$ where $P=L\cup \bigcup_{x\in H(L)}\Delta(x)$.}
\end{definition}

Now, let $\mathfrak{E}(L)$ be the set of all the finite extending standard forms of $L$.
In Fig.5, one can check that $P_2, P_3 \in \mathfrak{E}(L)$ but $P_1 \notin \mathfrak{E}(L)$.

In what follows, we write $L\hookrightarrow^{\prec} P$ when $L$ is a cover-preserving sublattice of a lattice $P$, and symbols $L\hookrightarrow P$,
$L\hookrightarrow^{\vee} P$ and $L\hookrightarrow^{\wedge} P$ stand for that $L$ is a sublattice, a $\vee$-subsemilattice and
a $\wedge$-subsemilattice of a lattice $P$, respectively. Then the following lemma is obvious.

\begin{lemma}\label{t31}
Suppose that $P\in \mathfrak{E}(L)$. Then $P$ is an
atomistic lattice, $\ell(L)=\ell(P)$ and $L\hookrightarrow^{\prec} P$.
\end{lemma}

For convenience, if $P\in \mathfrak{E}(L)$, then we denote $\Delta_P(L)=P-L$.
It is well known that a finite semimodular lattice $L$ can also
be expressed as sets of set (see \cite{Buchi52}). Therefore, by Lemma \ref{t31} and Definition \ref{d31}, if $P\in \mathfrak{E}(L)$, then there exists a lattice $(\mathcal{T}^{P}_L, \subseteq)$ with $\mathcal{T}^{P}_L\subseteq \mathcal{S}_{P}$ such that
\begin{equation}\label{e31}
L\cong \mathcal{T}^{P}_L \hookrightarrow^{\prec }\mathcal{S}_P\cong P
 \end{equation}
and
\begin{equation}\label{wwwa}\mbox{the identity map }i_d\mbox{ is a cover-preserving embedding map from }\mathcal{T}^{P}_{L}\mbox{ to }\mathcal{S}_{P}.\end{equation}
In fact, $\mathcal{T}^{P}_{L}\cup A(\mathcal{S}_{P})=\mathcal{S}_{P}$.

Consider the semimodular lattice $L$ and $L$'s extending standard form $P_2$ represented in Fig.5 again.
Then the two lattices $(\mathcal{T}^{P_2}_L, \subseteq)$ and $(\mathcal{S}_{P_2}, \subseteq)$ in Fig.6 satisfy
formula (\ref{e31}).

\par\noindent\vskip50pt
\begin{minipage}{11pc}
\setlength{\unitlength}{0.75pt}\begin{picture}(600,210)
\put(150,40){\circle{6}}\put(144,28){\makebox(0,0)[l]
{\footnotesize $\emptyset$}}
\put(150,80){\circle{6}}\put(125,78){\makebox(0,0)[l]
{\footnotesize$\{a\}$}}
\put(150,120){\circle{6}}\put(155,118){\makebox(0,0)[l]
{\footnotesize $\{b^{'},a\}$}}
\put(150,200){\circle{6}}\put(115,213){\makebox(0,0)[l]
{\footnotesize $\{d^{'},b^{'},a,c^{'}\}$}}
\put(190,160){\circle{6}}\put(195,168){\makebox(0,0)[l]
{\footnotesize $\{b^{'},a,c^{'}\}$}}
\put(110,160){\circle{6}}\put(65,175){\makebox(0,0)[l]
{\footnotesize $\{d^{'},b^{'},a\}$}}

\put(150,43){\line(0,1){34}}
\put(150,83){\line(0,1){34}}
\put(152,122){\line(1,1){36}}
\put(148,122){\line(-1,1){36}}
\put(112,162){\line(1,1){36}}
\put(188,162){\line(-1,1){36}}
\put(134,0){$(\mathcal{T}^{P_2}_L,\subseteq)$}

\put(370,40){\circle{6}}\put(364,28){\makebox(0,0)[l]
{\footnotesize $\emptyset$}}
\put(370,80){\circle{6}}\put(375,78){\makebox(0,0)[l]
{\footnotesize$\{a\}$}}
\put(370,120){\circle{6}}\put(328,120){\makebox(0,0)[l]
{\footnotesize $\{b^{'},a\}$}}
\put(370,200){\circle{6}}\put(335,215){\makebox(0,0)[l]
{\footnotesize $\{d^{'},b^{'},a,c^{'}\}$}}
\put(410,160){\circle{6}}\put(415,158){\makebox(0,0)[l]
{\footnotesize $\{b^{'},a,c^{'}\}$}}
\put(330,160){\circle{6}}\put(283,175){\makebox(0,0)[l]
{\footnotesize $\{d^{'},b^{'},a\}$}}
\put(330,80){\circle{6}}\put(305,88){\makebox(0,0)[l]
{\footnotesize $\{b^{'}\}$}}
\put(326,76){$\bullet$}
\put(290,80){\circle{6}}\put(260,88){\makebox(0,0)[l]
{\footnotesize $\{d^{'}\}$}}
\put(286,76){$\bullet$}
\put(410,80){\circle{6}}\put(415,78){\makebox(0,0)[l]
{\footnotesize $\{c^{'}\}$}}
\put(406,76){$\bullet$}
\put(370,43){\line(0,1){34}}
\put(370,83){\line(0,1){34}}
\put(372,122){\line(1,1){36}}
\put(368,122){\line(-1,1){36}}
\put(368,122){\line(-1,1){36}}
\put(332,162){\line(1,1){36}}
\put(408,162){\line(-1,1){36}}
\put(368,42){\line(-1,1){36}}
\put(332,82){\line(1,1){36}}
\put(368,42){\line(-2,1){74.5}}
\put(292,82){\line(1,2){37.5}}
\put(372,42){\line(1,1){36}}
\put(410,83){\line(0,1){74}}
\put(336,0){$(\mathcal{S}_{P_2}, \subseteq)$}
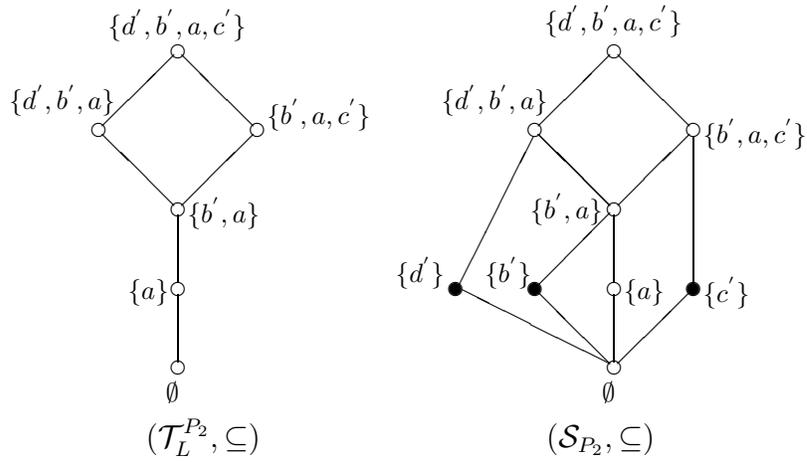
\put(115,-30){ Fig.6 Two lattices $(\mathcal{T}^{P_2}_L, \subseteq)$ and $(\mathcal{S}_{P_2}, \subseteq)$.}
\end{picture}
\end{minipage}
$$\quad$$
\\

By formula (\ref{e31}) and the construction of $L$'s extending standard forms, the following lemma is clearly.

\begin{lemma}\label{lemma}
\emph{Let $Q\in \mathfrak{E}(L)$ with $|A(Q)|>|J(L)|$. Then there exists an element $r\in \Delta_Q(L)$ such that $Q-\{r\}\in \mathfrak{E}(L)$. Further, let $P=Q-\{r\}$. Then we have that $L\cong \mathcal{T}^{P}_{L}\hookrightarrow^{\prec} \mathcal{S}_{P}$ and
$\mathcal{T}^{P}_{L}=\{X-\{r\}: X\in \mathcal{T}^{Q}_L\}$ where $\mathcal{T}^{Q}_L$ satisfies formulas (\ref{e31}) and (\ref{wwwa}).}
\end{lemma}

 The following example illustrates Lemma \ref{lemma}.
\begin{example}\label{exam1}
\emph{Consider the semimodular lattice $L$ and $L$'s extending standard form $P_3$ represented in Fig.5 again. Then the two lattices $(\mathcal{T}^{P_3}_L, \subseteq)$ and $(\mathcal{S}_{P_3}, \subseteq)$ in Fig.7 satisfy
formulas (\ref{e31}) and (\ref{wwwa}).}

\par\noindent\vskip50pt
\begin{minipage}{11pc}
\setlength{\unitlength}{0.75pt}\begin{picture}(600,170)
\put(110,40){\circle{6}}\put(105,28){\makebox(0,0)[l]
{\footnotesize $\emptyset$}}
\put(110,80){\circle{6}}\put(85,78){\makebox(0,0)[l]
{\footnotesize$\{a\}$}}
\put(110,120){\circle{6}}\put(115,118){\makebox(0,0)[l]
{\footnotesize $\{b^{'},a\}$}}
\put(110,200){\circle{6}}\put(65,215){\makebox(0,0)[l]
{\footnotesize $\{d^{'},d_{1}^{'},b^{'},a,c^{'}\}$}}
\put(150,160){\circle{6}}\put(155,168){\makebox(0,0)[l]
{\footnotesize $\{b^{'},a,c^{'}\}$}}
\put(70,160){\circle{6}}\put(5,175){\makebox(0,0)[l]
{\footnotesize $\{d^{'},d_{1}^{'},b^{'},a\}$}}

\put(110,43){\line(0,1){34}}
\put(110,83){\line(0,1){34}}
\put(112,122){\line(1,1){36}}
\put(108,122){\line(-1,1){36}}
\put(72,162){\line(1,1){36}}
\put(148,162){\line(-1,1){36}}
\put(85,0){$(\mathcal{T}^{P_3}_L,\subseteq)$}

\put(400,40){\circle{6}}\put(399,26){\makebox(0,0)[l]
{\footnotesize $\emptyset$}}
\put(400,80){\circle{6}}\put(405,78){\makebox(0,0)[l]
{\footnotesize$\{a\}$}}
\put(400,120){\circle{6}}\put(358,120){\makebox(0,0)[l]
{\footnotesize $\{b^{'},a\}$}}
\put(400,200){\circle{6}}\put(360,215){\makebox(0,0)[l]
{\footnotesize $\{d^{'},d_{1}^{'},b^{'},a,c^{'}\}$}}
\put(440,160){\circle{6}}\put(445,158){\makebox(0,0)[l]
{\footnotesize $\{b^{'},a,c^{'}\}$}}
\put(360,160){\circle{6}}\put(295,175){\makebox(0,0)[l]
{\footnotesize $\{d^{'},d_{1}^{'},b^{'},a\}$}}
\put(356,76){$\bullet$}
\put(360,80){\circle{6}}\put(337,90){\makebox(0,0)[l]
{\footnotesize $\{b^{'}\}$}}
\put(316,76){$\bullet$}
\put(320,80){\circle{6}}\put(295,90){\makebox(0,0)[l]
{\footnotesize $\{d^{'}\}$}}
\put(236,76){$\bullet$}
\put(240,80){\circle{6}}\put(210,88){\makebox(0,0)[l]
{\footnotesize $\{d_{1}^{'}\}$}}
\put(436,76){$\bullet$}
\put(440,80){\circle{6}}\put(445,78){\makebox(0,0)[l]
{\footnotesize $\{c^{'}\}$}}
\put(400,43){\line(0,1){34}}
\put(400,83){\line(0,1){34}}
\put(402,122){\line(1,1){36}}
\put(398,122){\line(-1,1){36}}
\put(398,122){\line(-1,1){36}}
\put(362,162){\line(1,1){36}}
\put(438,162){\line(-1,1){36}}
\put(398,42){\line(-1,1){36}}
\put(398,42){\line(-4,1){155}}
\put(242,82){\line(3,2){115}}
\put(362,82){\line(1,1){36}}
\put(398,42){\line(-2,1){74.5}}
\put(322,82){\line(1,2){37.5}}
\put(402,42){\line(1,1){36}}
\put(440,83){\line(0,1){74}}
\put(374,0){$(\mathcal{S}_{P_3}, \subseteq)$}
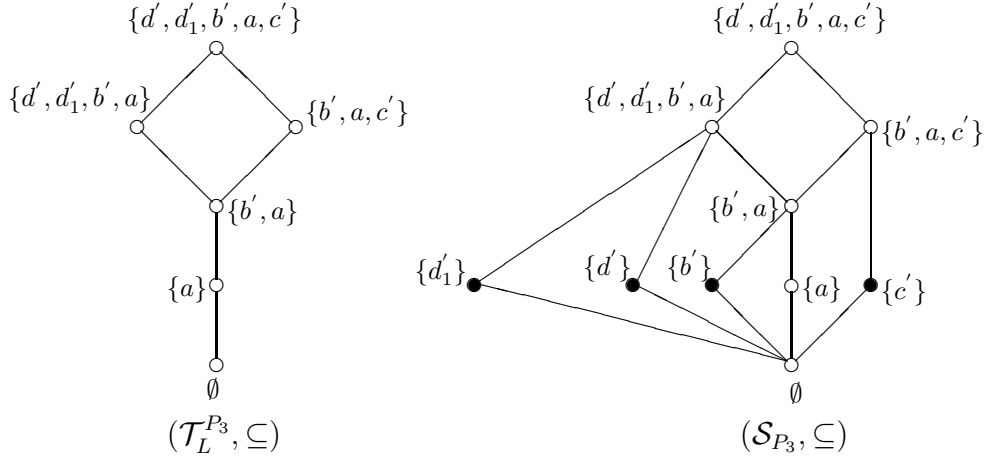
\put(110,-30){ \emph{Fig.7 Two lattices $(\mathcal{T}^{P_3}_L, \subseteq)$ and $(\mathcal{S}_{P_3}, \subseteq)$}.}
\end{picture}
\end{minipage}
$$\quad$$
\\

\emph{Obviously, $P_2=P_3 -\{d_{1}^{'}\}$. Then, from Fig.6, we know that $L\cong \mathcal{T}^{P_2}_{L}\hookrightarrow^{\prec} \mathcal{S}_{P_2}$ and $\mathcal{T}^{P_2}_{L}=\{X-\{d_{1}^{'}\}: X\in \mathcal{T}^{P_3}_L\}$.}
\end{example}

As a conclusion of this section, we shall supply an algorithm to construct a finite geometric lattice
$G$ which satisfies that $L\hookrightarrow^{\prec} G$ and $\ell(G)=\ell(L)$.

In the following, for each finite atomistic partially ordered set $P$ with $\ell(P)=m$,
we define two maps $\phi_P$ and $\varphi_P$ from $\{1,\cdots, m\}$ to the power set of $\mathcal{P}(A(P))$ and $\mathcal{P}(P)$ as
\begin{equation*}
\phi_P(i)=\{\sigma\in \mathfrak{I}_P(x): \ell(x)=i, x\in P\}
\end{equation*}
and
\begin{equation*}
\varphi_P(i)=\{x\in P: \ell(x)=i\},
\end{equation*}
respectively. Let $(\mathcal{O}, \subseteq)$ be an atomistic partially ordered set, and
$$\overline {X_\mathcal{O}}=\bigvee_{\mathcal{O}}\{T\in A(\mathcal{O}): T\subseteq X\}$$
for any set $X\subseteq 1_{\mathcal{O}}$ when $\bigvee_{\mathcal{O}}\{T\in A(\mathcal{O}): T\subseteq X\}$ exists.
Clearly, if $X\in \mathcal{O}$, then $\overline {X_\mathcal{O}}=X$.

Suppose that $P\in \mathfrak{E}(L)$ and $\ell(L)=m$. Then the following algorithm's output is a finite geometric lattice whose proof will be given in the next section.
\begin{algorithm}\label{a1}
\textbf{Input:} $\mathcal{Q}=\emptyset, \mathcal{R}=\mathcal{S}_P$, $k=3$, $t=0$ and $m=\ell(L)$. \textbf{Output:} $\mathcal{Q}$.\\
\emph{Step 1}. $\mathcal{Q}:=\mathcal{R}$ and $t:=k$. If there exists $X\in \varphi_{\mathcal{Q}}(k)$
which has a proper subset $U$ satisfying the following three conditions:

\emph{(i1)} if $V \in \varphi_{\mathcal{Q}}(t-1)$ and $Y\subseteq U\cap V$ then
$\overline {Y_\mathcal{R}}\subseteq U$;

\emph{(i2)} if $\sigma \in \phi_{\mathcal{Q}}(k-1)$ then $\bigcup\sigma\nsubseteq U$; and

\emph{(i3)} if $V \in \varphi_{\mathcal{Q}}(k-1)$, then $U\nsubseteq V$;\\
then $\mathcal{Q}:=\mathcal{Q}\cup \{U\}$. Otherwise, $k:=k+1$, and if $k\geq m+1$ then go to Step 5 and if not, go to Step 1.\\
\emph{Step 2}. If $\ell_{\mathcal{Q}}(U)=k-1$, then $k:=3$, $\mathcal{R}:=\mathcal{Q}$ and
go to Step 1. \\
\emph{Step 3}. If $U$ has a proper subset $W$ which satisfies the following three conditions:

\emph{(j1)} if $V \in \varphi_{\mathcal{Q}}(t-1)$ and $Y\subseteq W\cap V$ then
$\overline {Y_\mathcal{R}} \subseteq W$;

\emph{(j2)} if $\sigma \in \phi_{\mathcal{Q}}(k-2)$ then $\bigcup\sigma\nsubseteq W$; and

\emph{(j3)} if $V \in \varphi_{\mathcal{Q}}(k-2)$ then $W\nsubseteq V$; \\
then $\mathcal{Q}:=\mathcal{Q}\cup \{W\}$. Otherwise, go to Step 1.\\
\emph{Step 4}. If $\ell_{\mathcal{Q}}(W)=k-2$, then $k:=3$, $\mathcal{R}:=\mathcal{Q}$ and
go to Step 1. Otherwise, $k:=k-1$ and go to Step 3.\\
\emph{Step 5}. Stop.
\end{algorithm}
\section{All the finite geometric lattices}
In this section, we shall first prove that the output $\mathcal{Q}$ in
Algorithm \ref{a1} is an atomistic lattice, and then verify that $L$ is a cover-preserving
sublattice of $\mathcal{Q}$. Finally, we shall show that all the extending cover-preserving geometric
lattices of $L$ with the same length can be constructed by Algorithm \ref{a1}.

Below this paper, for convenience, if $(P, \leq)$ is a finite atomistic lattice with $n$ atoms, then we denote
$A(P)=\{1, \cdots, n\}$, and if $(\mathcal{S}, \subseteq)$ is a finite atomistic lattice with $m$ atoms,
then we denote $A(\mathcal{S})=\{\{1\},  \cdots, \{m\}\}$, and observe that
$$U\wedge_{\mathcal{S}} V=U\cap V$$ for any $U, V\in \mathcal{S}$.
\begin{lemma}\label{t32}
Every output $\mathcal{Q}$ in Algorithm \ref{a1} is a finite atomistic lattice.
\end{lemma}
\proof
Note that inasmuch as Algorithm \ref{a1} and Lemmas \ref{t21} and \ref{t31}, we know that every output $\mathcal{Q}$ is a finite
atomistic partially ordered set and it has the minimum element and the maximum element.
Then it suffices to show that the output $\mathcal{Q}$ is a $\wedge$-semilattice.
Note that $\mathcal{S}_P$ is an atomistic lattice by Lemma \ref{t31}. Obviously, the $\mathcal{R}$ in
Step 3 equals to the $\mathcal{R}$ in Step 1. Hence, by Algorithm \ref{a1},
we only need to prove that each partially ordered set $\mathcal{R}$ from Steps 2 and 4
returning to Step 1 in Algorithm \ref{a1} is a $\wedge$-semilattice.

For convenience, we next denote
$$\mathcal{D}_{H, G}=\{W\in \mathcal{R}: W\subseteq H \cap G\}.$$
One can see that $\mathcal{D}_{H, G}=\mathcal{D}_{G,H}$.

The rest of the proof will be completed in three steps.

A. If $\mathcal{R}$ is in Step 2, then $\mathcal{R}=\mathcal{Q}\cup \{M_{t-1}\}$,
in which $\mathcal{Q}$ is an atomistic lattice and $M_{t-1}$ is a proper subset of a certain $X$ in $\varphi_{\mathcal{Q}}(t)$ where $M_{t-1}$ satisfies the conditions (i1), (i2) and (i3). Let $E,F\in \mathcal{R}$.
Then there are three cases.

Case 1. If $E, F\in \mathcal{Q}$, then as $\mathcal{Q}$ is an atomistic lattice,
we know that $E\wedge_{\mathcal{Q}}F =E\cap F$ is the maximum of $\mathcal{D}_{E,F}$.
Therefore, $E\wedge_{\mathcal{R}} F=E\cap F\in \mathcal{R}$.

Case 2. If $E\in \mathcal{Q}$ and $F=M_{t-1}$, then suppose that $E\nparallel M_{t-1}$. Thus $E\wedge_{\mathcal{R}} M_{t-1}=E\cap M_{t-1}\in \mathcal{R}$ clearly. Now, assume that $E\parallel M_{t-1}$.
If $R\in \mathcal{D}_{E, M_{t-1}}$ then $R\subsetneq E$ and $R\subsetneq M_{t-1}$. We claim that
$X\nsubseteq E$. Otherwise $X\subseteq E$, which means that $M_{t-1}\subsetneq X\subseteq E$,
contrary to $E\parallel M_{t-1}$. We can distinguish two subcases.

Subcase 1. If $E\subsetneq X$, then $\ell_{\mathcal{Q}}(E)\leq t-1$ since
$X\in \varphi_{\mathcal{Q}}(t)$.

Subcase 2. If $E\parallel X$, then $E\wedge_{\mathcal{Q}} X=E\cap X\subsetneq X$.
Thus $\ell_{\mathcal{Q}}(E\cap X)\leq t-1$.

From Subcases 1 and 2, we know that there exists an element $E\cap X\in \mathcal{Q}$ such that
$\ell_{\mathcal{Q}}(E\cap X)\leq t-1$
and $E\cap M_{t-1}\subseteq E\cap X$ since $M_{t-1}\subsetneq X$.
Thus $$E\cap M_{t-1}\subseteq M_{t-1}\cap E\cap X,$$ and there exists an element $K\in \varphi_\mathcal{Q}(t-1)$ such that
$E\cap X\subseteq K$, or $E\cap X\in A(\mathcal{Q})$ by Algorithm \ref{a1}.
Then by (i1) of Algorithm \ref{a1}, we have that
$\overline{(E\cap M_{t-1})_{\mathcal{Q}}}\subseteq M_{t-1}$. Note that
$\overline{(E\cap M_{t-1})_{\mathcal{Q}}}\subseteq E$.
Therefore, $$\overline{(E\cap M_{t-1})_{\mathcal{Q}}}=E\cap M_{t-1}
\in \mathcal{D}_{E,M_{t-1}}.$$ Consequently,
$E\wedge_{\mathcal{R}} M_{t-1}=E\cap M_{t-1} \in \mathcal{R}$, i.e., $E\wedge_{\mathcal{R}} F \in \mathcal{R}$.

Case 3. If $E, F\in \mathcal{R}-\mathcal{Q}$, then clearly $E \wedge_{\mathcal{R}} F=E\cap F\in \mathcal{R}$.

In summary, $\mathcal{R}$ is a finite $\wedge$-semilattice.

B. If $\mathcal{R}$ is in Step 4 and $\mathcal{R}=\mathcal{Q}
\cup \{M_{t-1}, M_{t-2}\}$ in which $M_{t-2}$ is a proper subset of $M_{t-1}$ and it satisfies the conditions (j1), (j2) and (j3). By Algorithm \ref{a1}, we know that
$M_{t-2}\subsetneq  M_{t-1}\subsetneq X$. Suppose that $E, F\in \mathcal{R}$.
There are four cases as follows.

Case i. If $E, F\in \mathcal{Q}$, then similar to the proof of Case 1, we have that
$E\wedge_{\mathcal{R}} F=E\cap F\in \mathcal{R}$.

Case ii. If $E\in \mathcal{Q}$ and $F=M_{t-1}$, then similar to the proof of Case 2,
we know that $$E\wedge_{\mathcal{R}} M_{t-1}=E\cap  M_{t-1} \in \mathcal{R}.$$

Case iii. If $E\in \mathcal{Q}$ and $F=M_{t-2}$, then suppose that $E\nparallel M_{t-2}$.
Thus $$ E\wedge_{\mathcal{R}} M_{t-2}=E\cap M_{t-2}\in \mathcal{R}.$$
Now, assume that $E\parallel M_{t-2}$. Obviously, $R\subsetneq E$ and $R\subsetneq M_{t-2}$ for any $R\in  \mathcal{D}_{E, M_{t-2}}$, and $E\cap M_{t-2}\subsetneq E$. There are two subcases as follows.

Subcase (i). If $E\nparallel M_{t-1}$, then $M_{t-1}\subsetneq E$
or $E\subsetneq M_{t-1}$ since $E\neq M_{t-1}$. We claim that
$E\subsetneq M_{t-1}$.
Otherwise, $M_{t-2} \subsetneq M_{t-1}\subsetneq E$, contrary to the
fact that $E\parallel M_{t-2}$. Thus $E\subsetneq X$, and it follows from $X\in \varphi_{\mathcal{Q}}(t)$ that
$$\ell_{\mathcal{Q}}(E)\leq t-1.$$

Subcase (ii). If $E\parallel M_{t-1}$, then similar to the proof of Subcase 2, we know that there exists
an element $E\cap X\in \mathcal{Q}$ such that $$M_{t-2}\cap E\subseteq E\cap X \mbox{ and }
\ell_{\mathcal{Q}}(E\cap X)\leq t-1$$
since $M_{t-2}\cap E \subseteq M_{t-1}\cap E$.

Subcases (i) and (ii) mean that there exists an element $E\cap X\in \mathcal{Q}$ such that
$$\ell_{\mathcal{Q}}(E\cap X)\leq t-1\mbox{ and }E\cap M_{t-2}\subseteq E\cap X.$$
Hence $$E\cap M_{t-2}\subseteq M_{t-2}\cap E\cap X,$$ and there exists an element $K\in \varphi_\mathcal{Q}(t-1)$ such that
$E\cap X\subseteq K$, or $E\cap X\in A(\mathcal{Q})$ by Algorithm \ref{a1}.
Then by (j1) of Algorithm \ref{a1}, we have that
$\overline{(E\cap M_{t-2})_{\mathcal{Q}}}\subseteq M_{t-2}$. Note that
$\overline{(E\cap M_{t-2})_{\mathcal{Q}}}\subseteq E$.
Therefore, $$\overline{(E\cap M_{t-2})_{\mathcal{Q}}}=E\cap M_{t-2}
\in \mathcal{D}_{E,M_{t-2}}.$$ Consequently,
$E\wedge_{\mathcal{R}} M_{t-2}=E\cap M_{t-2} \in \mathcal{R}$.

Case iv. If $E, F\in \mathcal{R}-\mathcal{Q}$, then, clearly,
$E \wedge_{\mathcal{R}} F=E\cap F\in \mathcal{R}$.

In summary, $\mathcal{R}$ is a finite $\wedge$-semilattice.

C. Analogously, if $\mathcal{R}$ is in Step 4 and $\mathcal{R}=\mathcal{Q}
\cup \{M_{t-1}, M_{t-2}, \cdots, M_{t-r}\}$ for $r\in \{3,\cdots, t-2\}$ where $M_{t-r}$ is a proper subset of $M_{t-(r-1)}$ and it satisfies the conditions (j1), (j2) and (j3), then we can prove that $\mathcal{R}$ is a finite $\wedge$-semilattice.

To sum up, the output $\mathcal{Q}$ in Algorithm \ref{a1} is a finite atomistic lattice.
This completes the proof.
\endproof

Algorithm \ref{a1}, Definition \ref{d23} and Lemma \ref{t32} imply the following lemma.

\begin{lemma}\label{l31}
Let $P\in \mathfrak{E}(L)$, and $\mathcal{Q}$ be the output of Algorithm \ref{a1}. Then the following three statements hold.\\
\emph{(1)} $\ell_{\mathcal{S}_P}(X)=\ell_{\mathcal{Q}}(X)$ for any $X\in \mathcal{S}_P$.\\
\emph{(2)} If $\sigma \in \mathfrak{I}_{\mathcal{S}_P} (X)$, then $\sigma \in \mathfrak{I}_{\mathcal{Q}}(X)$
for any $X\in \mathcal{S}_P$.\\
\emph{(3)} If $X, Y\in \varphi_{\mathcal{Q}}(k)$ and $X\neq Y$, then
$\bigcup \sigma\nsubseteq Y$ for any $\sigma\in \mathfrak{I}_{\mathcal{Q}}(X)$.
\end{lemma}
\begin{lemma}\label{t33}
For every output $\mathcal{Q}$ in Algorithm \ref{a1}, $L\hookrightarrow^{\prec} \mathcal{Q}$.
\end{lemma}
\proof
Note that $\ell(L)=\ell(P)$ since $P\in \mathfrak{E}(L)$, and by Lemma \ref{l31},
$\ell(\mathcal{Q})=\ell(\mathcal{S}_{P})$.
Thus we have that $$\ell(\mathcal{Q})=\ell(L) \mbox{ since } P\cong \mathcal{S}_{P}.$$
By formula (\ref{e31}) and Algorithm \ref{a1}, there exists a lattice $\mathcal{T}^{P}_L$
such that $\mathcal{T}^{P}_L\subseteq \mathcal{S}_P \subseteq \mathcal{Q}$
and $L\cong \mathcal{T}^{P}_L \hookrightarrow^{\prec }\mathcal{S}_P$.
Therefore, we only need to prove that $\mathcal{T}^{P}_L \hookrightarrow\mathcal{Q}$.

First, by Lemma \ref{t31},
$\mathcal{S}_P$ is a finite atomistic lattice. Then
$E\wedge_{\mathcal{T}^{P}_{L}} F =E\cap F$ for any $E, F\in \mathcal{T}^{P}_{L}$ by (\ref{wwwa}).
On the other hand, from Lemma \ref{t32}, we know that $\mathcal{Q}$ is a finite atomistic lattice, which follows that
$H\wedge_{\mathcal{Q}} G=H\cap G$ for any $H, G\in \mathcal{Q}$.
Consequently, $\mathcal{T}^{P}_L \hookrightarrow^{\wedge }\mathcal{Q}$ since $\mathcal{T}^{P}_L\subseteq \mathcal{Q}$.

Next, we shall prove that $\mathcal{T}^{P}_L \hookrightarrow^{\vee }\mathcal{Q}$.

Let $M, N\in \mathcal{T}^{P}_L$.
Set $M\vee_{\mathcal{T}^{P}_L} N=Z$ and $M\vee_{\mathcal{Q}} N=T$.  Then $M\vee_{\mathcal{S}_P} N=Z$ by (\ref{wwwa}). Suppose that $T\neq Z$, then
$T\in \mathcal{Q}-\mathcal{S}_P$.
Obviously, $T\subsetneq Z$ since $\mathcal{T}^{P}_L\subseteq\mathcal{Q}$.  We claim that $M\parallel N$.
Otherwise, $T=Z=M\cup N$, a contradiction. As $\mathcal{T}^{P}_L$ is a finite semimodular lattice, we know that $\mathcal{T}^{P}_L$
contains a sublattice lattice as presented in Fig.8 (the required coverings $\prec$ and $\subseteq$ in the lattice $\mathcal{T}^{P}_L$ are
indicated by one line and double lines in Fig.8, respectively). Furthermore, by formula (\ref{e31}), Fig.8 is also a sublattice of $\mathcal{S}_P$ and
$\ell_{\mathcal{T}^{P}_L}(R)=\ell_{\mathcal{S}_P}(R)$ for any $R\in \mathcal{T}^{P}_L$. Therefore, by Lemma \ref{l31},
\begin{equation}\label{e34}
\ell_{\mathcal{Q}}(R)=\ell_{\mathcal{S}_P}(R)=\ell_{\mathcal{T}^{P}_L}(R)
\end{equation}
for
every $R\in \mathcal{T}^{P}_L$.
\par\noindent\vskip50pt
\begin{minipage}{11pc}
\setlength{\unitlength}{0.75pt}\begin{picture}(600,200)
\put(250,40){\circle{6}}\put(230,28){\makebox(0,0)[l]
{\footnotesize $M\cap N$}}
\put(250,80){\circle{6}}\put(223,78){\makebox(0,0)[l]
{\footnotesize$M_1$}}
\put(250,90){\makebox(0,0)[c]{$\cdot$}}
\put(250,100){\makebox(0,0)[c]{$\cdot$}}
\put(250,110){\makebox(0,0)[c]{$\cdot$}}
\put(250,120){\circle{6}}\put(223,122){\makebox(0,0)[l]
{\footnotesize $M_k$}}
\put(250,160){\circle{6}}\put(223,162){\makebox(0,0)[l]
{\footnotesize $M$}}

\put(290,80){\circle{6}}\put(303,78){\makebox(0,0)[l]
{\footnotesize$N$}}
\put(290,120){\circle{6}}\put(303,122){\makebox(0,0)[l]
{\footnotesize$N_1$}}
\put(290,130){\makebox(0,0)[c]{$\cdot$}}
\put(290,140){\makebox(0,0)[c]{$\cdot$}}
\put(290,150){\makebox(0,0)[c]{$\cdot$}}
\put(290,160){\circle{6}}\put(303,162){\makebox(0,0)[l]
{\footnotesize $N_k$}}
\put(290,200){\circle{6}}\put(285,212){\makebox(0,0)[l]
{\footnotesize $Z$}}

\put(250,43){\line(0,1){34}}
\put(250,123){\line(0,1){34}}
\put(290,163){\line(0,1){34}}
\put(290,83){\line(0,1){34}}
\put(252.5,41){\line(1,1){35.5}}
\put(252,43){\line(1,1){36}}
\put(252.5,81){\line(1,1){35.5}}
\put(252,83){\line(1,1){36}}
\put(252.5,121){\line(1,1){35.5}}
\put(252,123){\line(1,1){36}}
\put(252.5,161){\line(1,1){35.5}}
\put(252,163){\line(1,1){36}}

 \put(160,5){ Fig.8 A sublattice of $\mathcal{T}^{P}_L$.}
\end{picture}
\end{minipage}
\\

Now, set $\ell_{\mathcal{T}^{P}_L}(N)=t$. As $\mathcal{T}^{P}_L$ is a finite semimodular lattice and Fig.8
is a sublattice of $\mathcal{T}^{P}_L$, we obviously have that
$\ell_{\mathcal{T}^{P}_L}(N_k)=t+k$
and $\ell_{\mathcal{T}^{P}_L}(Z)=t+k+1$, which together with formula (\ref{e34}) imply
\begin{equation}\label{e35}\ell_{\mathcal{Q}}(N_k)=\ell_{\mathcal{S}_P}(N_k)=t+k \mbox{ and }\ell_{\mathcal{Q}}(Z)=\ell_{\mathcal{S}_P}(Z)=t+k+1.\end{equation}

Let $\eta \in \mathfrak{I}_{\mathcal{S}_{P}}(N)$.
Then by Definition \ref{d23} and formula (\ref{e35}), there exists a
subset $\rho$ of $A_{\mathcal{S}_{P}}(M)$ such that
$\eta \cup \rho=\pi\in \mathfrak{I}_{\mathcal{S}_{P}}(N_k)$
since Fig.8 is also a sublattice of $\mathcal{S}_{P}$.
Hence $\pi \in \mathfrak{I}_{\mathcal{Q}}(N_k)$ by Lemma \ref{l31}.

Using formula (\ref{e35}), $T\subsetneq Z$ and $T\in \mathcal{Q}-\mathcal{S}_P$, clearly,
there is a $T_0\in \mathcal{Q}-\mathcal{S}_P$ such that $M\cup N\subseteq T\subseteq T_0\subsetneq Z$
and $\ell_{\mathcal{Q}}(T_0)=t+k$, which follow by Lemma \ref{l31} that $\bigcup\sigma \nsubseteq T_0$ for every $\sigma\in \mathfrak{I}_{\mathcal{Q}}(N_k)$ since $\ell_{\mathcal{Q}}(N_k)=t+k$ and $T_0\neq N_k$. However,
$\pi=\eta\cup \rho\subseteq A_{\mathcal{S}_{P}}(M)\cup A_{\mathcal{S}_{P}}(N)$,
and then $\bigcup\pi \subseteq M\cup N\subseteq T_0$, contrary to $\pi \in \mathfrak{I}_{\mathcal{Q}}(N_k)$.

In summary, $\mathcal{T}^{P}_L \hookrightarrow^{\vee }\mathcal{Q}$.
Therefore, $L\hookrightarrow^{\prec} \mathcal{Q}$. This completes the proof.
\endproof

Notice that
\begin{equation}\label{wwwb}
\mbox{the identity map } i_d \mbox{ is a cover-preserving embedding map
from } \mathcal{T}^{P}_L \mbox{ to } \mathcal{Q}
\end{equation}
by the proof of Lemma \ref{t33}

Below, denote $$\mathfrak{S}=\{\mathcal{Q}: \mathcal{Q}
\mbox{ is an output of Algorithm \ref{a1}}\}$$ and
$$\overline{\mathfrak{S}}=\{\mathcal{Q}\in \mathfrak{S}: \mathcal{Q} \mbox{ satisfies the condition } (\mathfrak{M})\}$$
in which the condition $(\mathfrak{M})$ is as follows.\\
$(\mathfrak{M})$: If $X\in \varphi_{\mathcal{Q}}(k)$, then $\ell_{\mathcal{Q}}(\overline{(X\cup R)_{\mathcal{Q}}})=k+1$ for any $R\in A(\mathcal{Q})-A(X)$.

\begin{lemma}\label{t34}
Every $\mathcal{Q}\in \overline{\mathfrak{S}}$ is a finite geometric lattice.
\end{lemma}
\proof
By Lemma \ref{t32}, we know that $\mathcal{Q}$ is a finite atomistic lattice.
Then $[\emptyset, M]$ for any $M\in {\mathcal{Q}}$ is a geometric lattice when $\ell_{\mathcal{Q}}(M)\leq 2$.
Now, suppose that $[\emptyset, M]$ is a
geometric lattice for every $M\in \mathcal{Q}$ with $\ell_{\mathcal{Q}}(M)\leq k$. By induction, we shall prove that $[\emptyset, M]$ is a geometric lattice
for every $M\in \mathcal{Q}$ with $\ell_{\mathcal{Q}}(M)= k+1$.

Assume that $[\emptyset, M]$ is not a semimodular lattice. Then there exist two
elements $G, H\in [\emptyset, M]$ such that
$G\succ G\cap H, H\succ G\cap H$ but
$G\vee H\nsucc G$ or $G\vee H\nsucc H$,
say, $G\vee H\nsucc G$.
Note that $G\vee H\in [\emptyset, M]$. We claim that $G\vee H=M$. Otherwise, $G\vee H\subsetneq M$. Hence,
$\ell_{\mathcal{Q}}(G\vee H)\leq k$.
Therefore, $[\emptyset, G\vee H]$ is a geometric lattice. This follows that $G\vee H\succ G$ since $G, H, G\cap H, G\vee H\in [\emptyset, G\vee H]$, a contradiction.
Consequently, $G\vee H=M$, and which yields that $[\emptyset, M]$
contains a sublattice as presented in Fig.9 (the required coverings $\prec$ and $\subseteq$ in the lattice $[\emptyset, M]$ are indicated by one line and double lines in Fig.9, respectively).

\par\noindent\vskip50pt
\begin{minipage}{11pc}
\setlength{\unitlength}{0.75pt}\begin{picture}(600,220)
\put(250,40){\circle{6}}\put(230,28){\makebox(0,0)[l]
{\footnotesize $G\cap H$}}
\put(210,80){\circle{6}}\put(193,78){\makebox(0,0)[l]
{\footnotesize $H$}}

\put(290,80){\circle{6}}\put(303,78){\makebox(0,0)[l]
{\footnotesize$G$}}
\put(290,120){\circle{6}}\put(303,122){\makebox(0,0)[l]
{\footnotesize$G_1$}}
\put(290,130){\makebox(0,0)[c]{$\cdot$}}
\put(290,140){\makebox(0,0)[c]{$\cdot$}}
\put(290,150){\makebox(0,0)[c]{$\cdot$}}
\put(290,160){\circle{6}}\put(303,162){\makebox(0,0)[l]
{\footnotesize $G_{m-1}$}}
\put(290,200){\circle{6}}\put(303,202){\makebox(0,0)[l]
{\footnotesize $G_m$}}
\put(250,240){\circle{6}}\put(242,252){\makebox(0,0)[l]
{\footnotesize $M$}}
\put(248,43){\line(-1,1){36}}
\put(288,203){\line(-1,1){36}}
\put(290,83){\line(0,1){34}}
\put(252,43){\line(1,1){36}}
\put(290,163){\line(0,1){34}}
\put(210,83){\line(1,4){38.5}}
\put(213,82){\line(1,4){39}}
\put(160,5){Fig.9 A sublattice of $[\emptyset, M]$.}
\end{picture}
\end{minipage}
$$\quad$$

Let $R\in A(H)-A(G)$. We claim that $R\nsubseteq G_{m-1}$. Otherwise, $G_{m-1}\cap H \supsetneq G\cap H$. Then
$H\supseteq G_{m-1}\cap H\supsetneq G\cap H$ together with $H\succ G\cap H$
yields that $G_{m-1} \supseteq H$.
Thus $G\vee H\subseteq G_{m-1}$, a contradiction to the fact $G\vee H=M$. Therefore, $R\nsubseteq G_{m-1}$, then $R\in A(\mathcal{Q})-A(G_{m-1})$.
Note that $\ell_{\mathcal{Q}}(G_{m-1})\leq k-1$ by the
structure of Fig.9 and $\ell_{\mathcal{Q}}(M)=k+1$. Then by ($\mathfrak{M}$), we know that
\begin{equation}\label{e306}
      \ell_{\mathcal{Q}}(\overline{(G_{m-1}\cup R)_{\mathcal{Q}}})\leq k.
\end{equation}
On the other hand, $H\supseteq \overline{((G\cap H)\cup R)_{\mathcal{Q}}} \supsetneq G\cap H$, it follows
from $H\succ G\cap H$ that $H=\overline{((G\cap H)\cup R)_{\mathcal{Q}}}$.
Thus $\overline{(G_{m-1}\cup H)_{\mathcal{Q}}}=\overline{(G_{m-1}\cup R)_{\mathcal{Q}}}$, and then $\ell_{\mathcal{Q}}(\overline{(G_{m-1}\cup H)_{\mathcal{Q}}})\leq k$
by formula (\ref{e306}). However, $\overline{(G_{m-1}\cup H)_{\mathcal{Q}}}\supseteq G\vee H=M$ and $\ell_{\mathcal{Q}}(M)=k+1$, a contradiction. Therefore, $[\emptyset, M]$ is a semimodular lattice.

Consequently, $[\emptyset, M]$ is a
finite geometric lattice as $\mathcal{Q}$ is
a finite atomistic lattice, and the proof of the lemma is complete.
\endproof

Notice that from Lemmas \ref{t32}, \ref{t33} and \ref{t34}, we know that every output $\mathcal{Q}$ in
Algorithm \ref{a1} with condition ($\mathfrak{M}$) is a geometric lattice and $L$ is a cover-preserving sublattice of $\mathcal{Q}$.

The following example will illustrate that every output $\mathcal{Q}$ in Algorithm \ref{a1} with condition ($\mathfrak{M}$) is a geometric lattice and $L$ is a cover-preserving sublattice of $\mathcal{Q}$.
\begin{example}\label{ex31}
\emph{Consider the lattices $L$ and $P\in \mathfrak{E}(L)$ represented in Fig.10, respectively.}
\end{example}

\par\noindent\vskip50pt
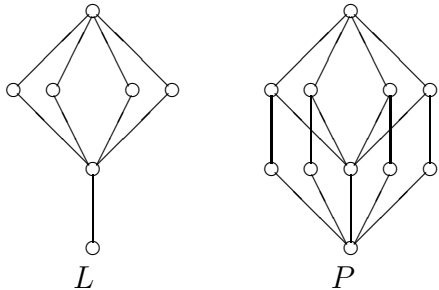
\begin{minipage}{11pc}
\setlength{\unitlength}{0.75pt}\begin{picture}(600,150)
\put(200,40){\circle{6}}
\put(200,80){\circle{6}}
\put(200,160){\circle{6}}
\put(240,120){\circle{6}}
\put(160,120){\circle{6}}
\put(220,120){\circle{6}}
\put(180,120){\circle{6}}

\put(200,43){\line(0,1){34}}
\put(202,82){\line(1,1){36}}
\put(198,82){\line(-1,1){36}}
\put(162,122){\line(1,1){36}}
\put(238,122){\line(-1,1){36}}
\put(202,82){\line(1,2){17.5}}
\put(198,82){\line(-1,2){17.5}}
\put(182,122){\line(1,2){17.5}}
\put(218,122){\line(-1,2){17.5}}
\put(190,20){$L$}

\put(330,40){\circle{6}}
\put(330,80){\circle{6}}
\put(330,160){\circle{6}}
\put(370,120){\circle{6}}
\put(290,120){\circle{6}}
\put(350,120){\circle{6}}
\put(310,120){\circle{6}}
\put(370,80){\circle{6}}
\put(290,80){\circle{6}}
\put(350,80){\circle{6}}
\put(310,80){\circle{6}}

\put(330,43){\line(0,1){34}}
\put(332,82){\line(1,1){36}}
\put(328,82){\line(-1,1){36}}
\put(292,122){\line(1,1){36}}
\put(368,122){\line(-1,1){36}}
\put(332,82){\line(1,2){17.5}}
\put(328,82){\line(-1,2){17.5}}
\put(312,122){\line(1,2){17.5}}
\put(348,122){\line(-1,2){17.5}}
\put(290,83){\line(0,1){34}}
\put(370,83){\line(0,1){34}}
\put(310,83){\line(0,1){34}}
\put(350,83){\line(0,1){34}}
\put(332,42){\line(1,2){17.5}}
\put(328,42){\line(-1,2){17.5}}
\put(332,42){\line(1,1){36}}
\put(328,42){\line(-1,1){36}}
\put(320,20){$P$}
\end{picture}
\end{minipage}

\par\noindent\vskip50pt
\begin{minipage}{11pc}
\setlength{\unitlength}{0.75pt}\begin{picture}(600,180)

\put(120,40){\circle{6}}\put(120,28){\makebox(0,0)[l]
{\footnotesize $\emptyset$}}
\put(120,100){\circle{6}}\put(126,98){\makebox(0,0)[l]
{\footnotesize $\{3\}$}}
\put(120,220){\circle{6}}\put(90,230){\makebox(0,0)[l]
{\footnotesize $\{1,2,3,4,5\}$}}
\put(220,160){\circle{6}}\put(225,158){\makebox(0,0)[l]
{\footnotesize $\{3,5\}$}}
\put(20,160){\circle{6}}\put(-18,158){\makebox(0,0)[l]
{\footnotesize $\{1,3\}$}}
\put(180,160){\circle{6}}\put(140,160){\makebox(0,0)[l]
{\footnotesize $\{3, 4\}$}}
\put(60,160){\circle{6}}\put(70,160){\makebox(0,0)[l]
{\footnotesize $\{2, 3\}$}}

\put(120,43){\line(0,1){54}}
\put(122,102){\line(1,1){56}}
\put(118,102){\line(-1,1){56}}

\put(62,162){\line(1,1){56}}
\put(178,162){\line(-1,1){56}}
\put(122,102){\line(5,3){95}}
\put(118,102){\line(-5,3){95}}
\put(22,162){\line(5,3){95}}
\put(218,162){\line(-5,3){95}}
\put(115,0){$\mathcal{T}^P_{L}$}

\put(400,40){\circle{6}}\put(395,28){\makebox(0,0)[l]
{\footnotesize $\emptyset$}}
\put(400,100){\circle{6}}\put(375,97){\makebox(0,0)[l]
{\footnotesize $\{3\}$}}
\put(400,220){\circle{6}}\put(370,230){\makebox(0,0)[l]
{\footnotesize $\{1,2,3,4,5\}$}}
\put(500,160){\circle{6}}\put(503,158){\makebox(0,0)[l]
{\footnotesize $\{3,5\}$}}
\put(500,100){\circle{6}}\put(503,98){\makebox(0,0)[l]
{\footnotesize $\{5\}$}}
\put(300,160){\circle{6}}\put(273,172){\makebox(0,0)[l]
{\footnotesize $\{1,3\}$}}
\put(300,100){\circle{6}}\put(274,98){\makebox(0,0)[l]
{\footnotesize $\{1\}$}}
\put(460,160){\circle{6}}\put(420,160){\makebox(0,0)[l]
{\footnotesize $\{3, 4\}$}}
\put(460,100){\circle{6}}\put(435,100){\makebox(0,0)[l]
{\footnotesize $\{4\}$}}
\put(340,160){\circle{6}}\put(351,160){\makebox(0,0)[l]
{\footnotesize $\{2, 3\}$}}
\put(340,100){\circle{6}}\put(343,102){\makebox(0,0)[l]
{\footnotesize $\{2\}$}}

\put(460,103){\line(0,1){54}}
\put(340,103){\line(0,1){54}}
\put(300,103){\line(0,1){54}}
\put(500,103){\line(0,1){54}}

\put(400,43){\line(0,1){54}}
\put(402,102){\line(1,1){56}}
\put(398,102){\line(-1,1){56}}

\put(342,162){\line(1,1){56}}
\put(458,162){\line(-1,1){56}}
\put(402,102){\line(5,3){95}}
\put(398,102){\line(-5,3){95}}
\put(302,162){\line(5,3){95}}
\put(498,162){\line(-5,3){95}}

\put(402,42){\line(5,3){95}}
\put(398,42){\line(-5,3){95}}
\put(402,42){\line(1,1){56}}
\put(398,42){\line(-1,1){56}}
\put(390,0){$\mathcal{S}_{P}$}

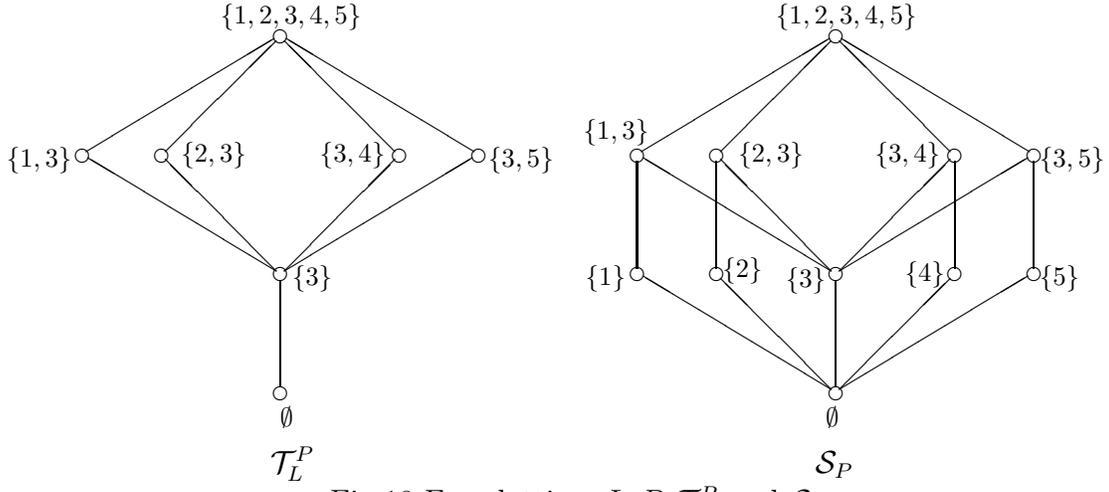
\put(140,-20){ Fig.10 Four lattices $L, P, \mathcal{T}^P_{L}$ and $\mathcal{S}_{P}$.}
\end{picture}
\end{minipage}
$$\quad$$
\\
Obviously, $\mathcal{T}^P_{L}$ and $\mathcal{S}_{P}$ satisfy formula (\ref{e31}), respectively.

$\textbf{Input}: \mathcal{Q}=\emptyset, \mathcal{R}=\mathcal{S}_P, k=3$, $t=0$ and $m=3$.

$\textbf{Output}: \mathcal{Q}.$

Step 1. $\mathcal{Q}:=\mathcal{R}$, $t:=3$, $U_1=\{1,2,4\}$ is a proper subset of $\{1,2,3,4,5\}$ satisfying (i1), (i2) and (i3), and $\mathcal{Q}:=\mathcal{Q}\cup \{U_1\}$.

Step 2. $\ell_{\mathcal{Q}}(U_1)=2$, $k:=3$ and $\mathcal{R}:=\mathcal{Q}$ (the lattice
$\mathcal{R}$ as represented in Fig.11).

Step 3. $\mathcal{Q}:=\mathcal{R}$, $t:=3$, $U_2=\{1,5\}$ is a proper subset of $\{1,2,3,4,5\}$ satisfying (i1), (i2) and (i3), and $\mathcal{Q}:=\mathcal{Q}\cup \{U_2\}$.

Step 4. $\ell_{\mathcal{Q}}(U_2)=2$, $k:=3$ and $\mathcal{R}:=\mathcal{Q}$ (the lattice
$\mathcal{R}$ as represented in Fig.12).

Step 5. $\mathcal{Q}:=\mathcal{R}$, $t:=3$, $U_3=\{2,5\}$ is a proper subset of $\{1,2,3,4,5\}$ satisfying (i1), (i2) and (i3), and $\mathcal{Q}:=\mathcal{Q}\cup \{U_3\}$.

Step 6. $\ell_{\mathcal{Q}}(U_3)=2$, $k:=3$ and $\mathcal{R}:=\mathcal{Q}$ (the lattice
$\mathcal{R}$ as represented in Fig.13).

Step 7. $\mathcal{Q}:=\mathcal{R}$, $t:=3$, $U_4=\{4,5\}$ is a proper subset of $\{1,2,3,4,5\}$ satisfying (i1), (i2) and (i3), and $\mathcal{Q}:=\mathcal{Q}\cup \{U_4\}$.

Step 8. $\ell_{\mathcal{Q}}(U_4)=2$, $k:=3$ and $\mathcal{R}:=\mathcal{Q}$ (the lattice
$\mathcal{R}$ as represented in Fig.14).

Step 9. $\mathcal{Q}:=\mathcal{R}$, $t:=3$ and $\{1,2,3,4,5\}$ has no proper subset satisfying (i1), (i2) and (i3), $k=4\geq 4$.

Step 10. Stop.

\par\noindent\vskip50pt
\begin{minipage}{11pc}
\setlength{\unitlength}{0.75pt}\begin{picture}(600,300)

\put(300,40){\circle{6}}\put(295,28){\makebox(0,0)[l]
{\footnotesize $\emptyset$}}
\put(300,120){\circle{6}}\put(275,118){\makebox(0,0)[l]
{\footnotesize $\{3\}$}}
\put(300,280){\circle{6}}\put(270,295){\makebox(0,0)[l]
{\footnotesize $\{1,2,3,4,5\}$}}
\put(460,200){\circle{6}}\put(463,198){\makebox(0,0)[l]
{\footnotesize $\{3,5\}$}}
\put(460,120){\circle{6}}\put(463,118){\makebox(0,0)[l]
{\footnotesize $\{5\}$}}
\put(140,200){\circle{6}}\put(104,198){\makebox(0,0)[l]
{\footnotesize $\{1,3\}$}}
\put(140,120){\circle{6}}\put(114,118){\makebox(0,0)[l]
{\footnotesize $\{1\}$}}
\put(380,200){\circle{6}}\put(385,200){\makebox(0,0)[l]
{\footnotesize $\{3, 4\}$}}
\put(380,120){\circle{6}}\put(385,120){\makebox(0,0)[l]
{\footnotesize $\{4\}$}}
\put(220,200){\circle{6}}\put(185,198){\makebox(0,0)[l]
{\footnotesize $\{2, 3\}$}}
\put(220,120){\circle{6}}\put(196,120){\makebox(0,0)[l]
{\footnotesize $\{2\}$}}

\put(300,200){\circle{6}}\put(304,200){\makebox(0,0)[l]
{\footnotesize $\{1,2,4\}$}}

\put(380,123){\line(0,1){74}}
\put(220,123){\line(0,1){74}}
\put(140,123){\line(0,1){74}}
\put(460,123){\line(0,1){74}}

\put(300,43){\line(0,1){74}}
\put(302,122){\line(1,1){76}}
\put(298,122){\line(-1,1){76}}

\put(222,202){\line(1,1){76}}
\put(378,202){\line(-1,1){76}}
\put(302,122){\line(2,1){155}}
\put(298,122){\line(-2,1){155}}
\put(142,202){\line(2,1){155}}
\put(458,202){\line(-2,1){155}}

\put(302,42){\line(2,1){155}}
\put(298,42){\line(-2,1){155}}
\put(302,42){\line(1,1){76}}
\put(298,42){\line(-1,1){76}}
\put(378,122){\line(-1,1){76}}
\put(222,122){\line(1,1){76}}
\put(142,122){\line(2,1){155}}
\put(300,203){\line(0,1){74}}

\put(290,5){$\mathcal{R}$}

\put(230,-20){ Fig.11 The lattice $\mathcal{R}$.}
\end{picture}
\end{minipage}
$$\quad$$
\\

\par\noindent\vskip50pt
\begin{minipage}{11pc}
\setlength{\unitlength}{0.75pt}\begin{picture}(600,290)

\put(300,40){\circle{6}}\put(295,28){\makebox(0,0)[l]
{\footnotesize $\emptyset$}}
\put(300,120){\circle{6}}\put(275,118){\makebox(0,0)[l]
{\footnotesize $\{3\}$}}
\put(300,280){\circle{6}}\put(270,295){\makebox(0,0)[l]
{\footnotesize $\{1,2,3,4,5\}$}}
\put(460,200){\circle{6}}\put(463,198){\makebox(0,0)[l]
{\footnotesize $\{3,5\}$}}
\put(460,120){\circle{6}}\put(463,118){\makebox(0,0)[l]
{\footnotesize $\{5\}$}}
\put(140,200){\circle{6}}\put(104,198){\makebox(0,0)[l]
{\footnotesize $\{1,3\}$}}
\put(140,120){\circle{6}}\put(114,118){\makebox(0,0)[l]
{\footnotesize $\{1\}$}}
\put(380,200){\circle{6}}\put(385,200){\makebox(0,0)[l]
{\footnotesize $\{3, 4\}$}}
\put(380,120){\circle{6}}\put(385,120){\makebox(0,0)[l]
{\footnotesize $\{4\}$}}
\put(220,200){\circle{6}}\put(185,198){\makebox(0,0)[l]
{\footnotesize $\{2, 3\}$}}
\put(220,120){\circle{6}}\put(196,120){\makebox(0,0)[l]
{\footnotesize $\{2\}$}}

\put(300,200){\circle{6}}\put(304,200){\makebox(0,0)[l]
{\footnotesize $\{1,2,4\}$}}

\put(60,200){\circle{6}}\put(24,200){\makebox(0,0)[l]
{\footnotesize $\{1,5\}$}}
\put(380,123){\line(0,1){74}}
\put(220,123){\line(0,1){74}}
\put(140,123){\line(0,1){74}}
\put(460,123){\line(0,1){74}}

\put(300,43){\line(0,1){74}}
\put(302,122){\line(1,1){76}}
\put(298,122){\line(-1,1){76}}

\put(222,202){\line(1,1){76}}
\put(378,202){\line(-1,1){76}}
\put(302,122){\line(2,1){155}}
\put(298,122){\line(-2,1){155}}
\put(142,202){\line(2,1){155}}
\put(458,202){\line(-2,1){155}}

\put(302,42){\line(2,1){155}}
\put(298,42){\line(-2,1){155}}
\put(302,42){\line(1,1){76}}
\put(298,42){\line(-1,1){76}}
\put(378,122){\line(-1,1){76}}
\put(222,122){\line(1,1){76}}
\put(142,122){\line(2,1){155}}
\put(300,203){\line(0,1){74}}
\put(62,202){\line(3,1){235}}
\put(290,5){$\mathcal{R}$}
\put(458,122){\line(-5,1){395}}
\put(138,122){\line(-1,1){76}}
\put(230,-20){ Fig.12 The lattice $\mathcal{R}$.}
\end{picture}
\end{minipage}
$$\quad$$
\\

\par\noindent\vskip50pt
\begin{minipage}{11pc}
\setlength{\unitlength}{0.75pt}\begin{picture}(600,300)

\put(300,40){\circle{6}}\put(295,28){\makebox(0,0)[l]
{\footnotesize $\emptyset$}}
\put(300,120){\circle{6}}\put(275,118){\makebox(0,0)[l]
{\footnotesize $\{3\}$}}
\put(300,280){\circle{6}}\put(270,295){\makebox(0,0)[l]
{\footnotesize $\{1,2,3,4,5\}$}}
\put(460,200){\circle{6}}\put(463,198){\makebox(0,0)[l]
{\footnotesize $\{3,5\}$}}
\put(460,120){\circle{6}}\put(463,118){\makebox(0,0)[l]
{\footnotesize $\{5\}$}}
\put(140,200){\circle{6}}\put(104,198){\makebox(0,0)[l]
{\footnotesize $\{1,3\}$}}
\put(140,120){\circle{6}}\put(114,118){\makebox(0,0)[l]
{\footnotesize $\{1\}$}}
\put(380,200){\circle{6}}\put(385,200){\makebox(0,0)[l]
{\footnotesize $\{3, 4\}$}}
\put(380,120){\circle{6}}\put(385,120){\makebox(0,0)[l]
{\footnotesize $\{4\}$}}
\put(220,200){\circle{6}}\put(185,198){\makebox(0,0)[l]
{\footnotesize $\{2, 3\}$}}
\put(220,120){\circle{6}}\put(196,120){\makebox(0,0)[l]
{\footnotesize $\{2\}$}}

\put(300,200){\circle{6}}\put(304,200){\makebox(0,0)[l]
{\footnotesize $\{1,2,4\}$}}

\put(60,200){\circle{6}}\put(24,200){\makebox(0,0)[l]
{\footnotesize $\{1,5\}$}}
\put(540,200){\circle{6}}\put(543,198){\makebox(0,0)[l]
{\footnotesize $\{2,5\}$}}

\put(380,123){\line(0,1){74}}
\put(220,123){\line(0,1){74}}
\put(140,123){\line(0,1){74}}
\put(460,123){\line(0,1){74}}

\put(300,43){\line(0,1){74}}
\put(302,122){\line(1,1){76}}
\put(298,122){\line(-1,1){76}}

\put(222,202){\line(1,1){76}}
\put(378,202){\line(-1,1){76}}
\put(302,122){\line(2,1){155}}
\put(298,122){\line(-2,1){155}}
\put(142,202){\line(2,1){155}}
\put(458,202){\line(-2,1){155}}

\put(302,42){\line(2,1){155}}
\put(298,42){\line(-2,1){155}}
\put(302,42){\line(1,1){76}}
\put(298,42){\line(-1,1){76}}
\put(378,122){\line(-1,1){76}}
\put(222,122){\line(1,1){76}}
\put(142,122){\line(2,1){155}}
\put(300,203){\line(0,1){74}}
\put(62,202){\line(3,1){235}}
\put(458,122){\line(-5,1){395}}
\put(138,122){\line(-1,1){76}}
\put(462,122){\line(1,1){76}}
\put(538,202){\line(-3,1){235}}
\put(222,122){\line(4,1){315}}

\put(290,5){$\mathcal{R}$}

\put(230,-20){ Fig.13 The lattice $\mathcal{R}$.}
\end{picture}
\end{minipage}
$$\quad$$
\\

\par\noindent\vskip50pt
\begin{minipage}{11pc}
\setlength{\unitlength}{0.75pt}\begin{picture}(600,300)

\put(300,40){\circle{6}}\put(295,28){\makebox(0,0)[l]
{\footnotesize $\emptyset$}}
\put(300,120){\circle{6}}\put(275,118){\makebox(0,0)[l]
{\footnotesize $\{3\}$}}
\put(300,280){\circle{6}}\put(270,295){\makebox(0,0)[l]
{\footnotesize $\{1,2,3,4,5\}$}}
\put(460,200){\circle{6}}\put(463,198){\makebox(0,0)[l]
{\footnotesize $\{3,5\}$}}
\put(460,120){\circle{6}}\put(463,118){\makebox(0,0)[l]
{\footnotesize $\{5\}$}}
\put(140,200){\circle{6}}\put(104,198){\makebox(0,0)[l]
{\footnotesize $\{1,3\}$}}
\put(140,120){\circle{6}}\put(114,118){\makebox(0,0)[l]
{\footnotesize $\{1\}$}}
\put(380,200){\circle{6}}\put(385,200){\makebox(0,0)[l]
{\footnotesize $\{3, 4\}$}}
\put(380,120){\circle{6}}\put(385,120){\makebox(0,0)[l]
{\footnotesize $\{4\}$}}
\put(220,200){\circle{6}}\put(185,198){\makebox(0,0)[l]
{\footnotesize $\{2, 3\}$}}
\put(220,120){\circle{6}}\put(196,120){\makebox(0,0)[l]
{\footnotesize $\{2\}$}}
\put(300,200){\circle{6}}\put(304,200){\makebox(0,0)[l]
{\footnotesize $\{1,2,4\}$}}
\put(60,200){\circle{6}}\put(24,200){\makebox(0,0)[l]
{\footnotesize $\{1,5\}$}}
\put(-20,200){\circle{6}}\put(-24,216){\makebox(0,0)[l]
{\footnotesize $\{4,5\}$}}
\put(540,200){\circle{6}}\put(543,198){\makebox(0,0)[l]
{\footnotesize $\{2,5\}$}}

\put(380,123){\line(0,1){74}}
\put(378,122){\line(-5,1){395}}
\put(458,122){\line(-6,1){475}}
\put(220,123){\line(0,1){74}}
\put(140,123){\line(0,1){74}}
\put(460,123){\line(0,1){74}}

\put(300,43){\line(0,1){74}}
\put(302,122){\line(1,1){76}}
\put(298,122){\line(-1,1){76}}

\put(222,202){\line(1,1){76}}
\put(378,202){\line(-1,1){76}}
\put(302,122){\line(2,1){155}}
\put(298,122){\line(-2,1){155}}
\put(142,202){\line(2,1){155}}
\put(458,202){\line(-2,1){155}}

\put(302,42){\line(2,1){155}}
\put(298,42){\line(-2,1){155}}
\put(302,42){\line(1,1){76}}
\put(298,42){\line(-1,1){76}}
\put(378,122){\line(-1,1){76}}
\put(222,122){\line(1,1){76}}
\put(142,122){\line(2,1){155}}
\put(300,203){\line(0,1){74}}
\put(62,202){\line(3,1){235}}
\put(458,122){\line(-5,1){395}}
\put(138,122){\line(-1,1){76}}
\put(462,122){\line(1,1){76}}
\put(538,202){\line(-3,1){235}}
\put(222,122){\line(4,1){315}}
\put(-18,202){\line(4,1){315}}
\put(290,5){$\mathcal{R}$}

\put(230,-20){ Fig.14 The lattice $\mathcal{R}$.}
\end{picture}
\end{minipage}
$$\quad$$
\\
Therefore, the output $\mathcal{Q}$ in Algorithm \ref{a1} is the $\mathcal{R}$ as represented in Fig.14.
One can check that $\mathcal{Q}\in \overline{\mathfrak{S}}$, $\mathcal{Q}$ is a finite geometric lattice and
$L\hookrightarrow^{\prec} \mathcal{Q}$.

\begin{definition}\label{d41}
\emph{Let $L_1$ and $L_2$ be two finite atomistic lattices with $L_1\subseteq L_2$.
If $L_1$ satisfies: for any $p\in L_1$, }\\
\emph{(e1)} $ A_{L_{1}}(p)=A_{L_{2}}(p)$;\\
\emph{(e2)} $\ell_{L_{1}}(p)=\ell_{L_{2}}(p)$; and\\
\emph{(e3)} $[0, p]_{L_1}=[0, p]_{L_{2}}
\emph{ when }\ell_{L_{1}}(p)\leq k$,\\
\emph{then we
say that $L_1$ is a $k$ order normal subset lattice of $L_2$.}
\end{definition}

\begin{lemma}\label{t35}
Let $(\mathcal{Q},\subseteq)$ be a finite geometric lattice with
$L\hookrightarrow^{\prec} \mathcal{Q}$ and $\ell(L)=\ell(\mathcal{Q})$.
Then $\mathcal{Q}\in \overline{\mathfrak{S}}$.
\end{lemma}
\proof
Let $\ell(\mathcal{Q})=m$. As $\mathcal{Q}$ is a finite geometric lattice, $\mathcal{Q}$ satisfies condition ($\mathfrak{M}$). Thus
we only need to prove that $\mathcal{Q}$ is an output of Algorithm \ref{a1}.
Since $L\hookrightarrow^{\prec} \mathcal{Q}$, there exists a lattice
$\mathcal{T}\subseteq \mathcal{Q}$ such that
$L\cong \mathcal{T}\hookrightarrow^{\prec} \mathcal{Q}$. Hence there exists a lattice $P\in \mathfrak{E}(L)$
such that $\mathcal{S}_{P}=\mathcal{T}\cup A(\mathcal{Q})$. Thus there exists a lattice $\mathcal{T}^{P}_L$
such that $\mathcal{T}^{P}_L\subseteq \mathcal{S}_P$ and
$L\cong \mathcal{T}^{P}_L \hookrightarrow^{\prec }\mathcal{S}_P$ by formula (\ref{e31}).

Because $\mathcal{Q}$ is geometric, the following four statements hold.\\
C1. For every $M\in \mathcal{Q}$, $[\emptyset, M]_{\mathcal{Q}}$ is a geometric lattice.\\
C2. If $M, N\in \mathcal{Q}$, then $\overline{(M\cap N)_{\mathcal{Q}}}=M\cap N$.\\
C3. If $M, N\in \varphi_{\mathcal{Q}}(k)$ and $M\neq N$,
then $M\nsubseteq N$ and $N\nsubseteq M$.\\
C4. If $\sigma\in \phi_{\mathcal{Q}}(k)$, $M\in \varphi_{\mathcal{Q}}(k)$ and
$\bigvee_{\mathcal{Q}}\sigma\neq M$, then $\sigma\notin \mathfrak{I}_{\mathcal{Q}}(M)$
and $\bigcup \sigma\nsubseteq M$.

The rest of the proof will be completed in three steps.

(I). Let $\mathcal{R}_{1}=\mathcal{S}_{P}$. Then by Definition \ref{d41}, $\mathcal{R}_{1}$ is a 2 order normal subset lattice of $\mathcal{Q}$.
Suppose that $X\in \mathcal{R}_{1}$ and $\ell_{\mathcal{R}_{1}}(X)=3$.
Let $U\in \mathcal{Q}-\mathcal{R}_{1}$, $U\subsetneq X$ and
$\ell_{\mathcal{Q}}(U)=2$. If $V \in \varphi_{\mathcal{R}_{1}}(2)$,
then $\overline{(U\cap V)_{\mathcal{R}_{1}}}=\overline{(U\cap V)_{\mathcal{Q}}}$
by (e3). It follows from C2 that
\begin{equation}\label{tian1}
\overline{(U\cap V)_{\mathcal{R}_{1}}}\subseteq U.
\end{equation}
Obviously, by C4 and (e2) and (e3) in Definition \ref{d41},
\begin{equation}\label{tian2}
\bigcup\sigma \nsubseteq U
\end{equation}
for any $\sigma\in \phi_{\mathcal{R}_{1}}(2)$. Moreover, by (e2) in Definition \ref{d41} and C3, we have
\begin{equation}\label{tian3}
U\nsubseteq V
\end{equation}
for every $V \in \varphi_{\mathcal{R}_{1}}(2)$ since $U\neq V$.
Thus, by formulas (\ref{tian1}), (\ref{tian2}) and (\ref{tian3}),
we know that $U$ satisfies (i1), (i2) and (i3) in Algorithm \ref{a1}. Therefore, $\mathcal{R}_{1}\cup \{U\}$ is an atomistic lattice by the proof of Lemma \ref{t32}. Clearly, $\mathcal{R}_{1}\cup \{U\}$ is a 2 order
normal subset lattice of $\mathcal{Q}$.

Suppose that $E\in \mathcal{R}_{1}\cup \{U\}$ and $\ell_{\mathcal{R}_{1}\cup \{U\}}(E)=3$.
Let $U_1\in \mathcal{Q}-(\mathcal{R}_{1}\cup \{U\})$, $U_1\subsetneq E$ and
$\ell_{\mathcal{Q}}(U_1)=2$.
Similar to the proof of the preceding paragraph, we can prove that $\mathcal{R}_{1}\cup \{U\}\cup \{U_1\}$ is an atomistic lattice which is a 2 order normal subset lattice of $\mathcal{Q}$ since
$\mathcal{R}_{1}\cup \{U\}$ is a 2 order normal subset lattice of $\mathcal{Q}$.

Repeating the process as above, we can obtain an atomistic lattice
$$\mathcal{R}_2=\mathcal{R}_1 \cup \bigcup\{[\emptyset, M]_\mathcal{Q}: \ell_{\mathcal{R}_1}(M)=3\}.$$
Obviously, \begin{equation}\label{equt}\mathcal{R}_2=\mathcal{R}_1 \cup \bigcup\{[\emptyset, M]_\mathcal{Q}: \ell_{\mathcal{Q}}(M)=3, M\in \mathcal{R}_1\}\end{equation} by (e2) in Definition \ref{d41}. Therefore, $\mathcal{R}_2$ is a 3 order normal subset lattice
of $\mathcal{Q}$, and for any $F\in \mathcal{R}_{2}$ with
$\ell_{\mathcal{R}_{2}}(F)\leq 3$, $F$ has no proper subset $N$ satisfying
(i1), (i2) and (i3) in Algorithm \ref{a1}.

(II). Suppose that $X\in \mathcal{R}_{2}$ and $\ell_{\mathcal{R}_{2}}(X)=4$. Let $U\in \mathcal{Q}-\mathcal{R}_2$, $U\subsetneq X$ and $\ell_{\mathcal{Q}}(U)=3$. There are two cases as below.

Case 1. If there exists $E\in \varphi_{\mathcal{R}_{2}}(2)$ such
that $E\subsetneq U$. Similar to the proof of formulas (\ref{tian1}), (\ref{tian2}) and (\ref{tian3}),
we can verify that $U$ satisfies (i1), (i2) and (i3) in Algorithm \ref{a1}. Therefore, $\mathcal{R}_{2}\cup \{U\}$ is an atomistic lattice by the proof of Lemma \ref{t32}. Clearly, $\mathcal{R}_{2}\cup \{U\}$ is a
2 order normal subset lattice of $\mathcal{Q}$. Thus, similar to the proof of (\ref{equt}), we can obtain an atomistic lattice $\mathcal{R}_{2}\cup [\emptyset, U]_{\mathcal{Q}}$ that is a 3 order normal subset lattice of $\mathcal{Q}$.

Case 2. If there is no element $E\in  \varphi_{\mathcal{R}_{2}}(2)$ such that $E\subsetneq U$,
then there exists $U_1\in \mathcal{Q}-\mathcal{R}_2$ such that $U_1\subsetneq U$ and
$\ell_{\mathcal{Q}}(U_1)=2$ since $\mathcal{Q}$ is geometric. Thus we have that the following three results.

(a1) By (e2), (e3), C1 and C2, $\overline{(U_1\cap V)_{\mathcal{R}_{2}}}\subseteq U_1$
for any $V\in \varphi_{\mathcal{R}_{2}}(3)$;

(a2) By (e2) and C4, $\bigcup \sigma\nsubseteq U_1$
for any $\sigma\in \phi_{\mathcal{R}_{2}}(2)$;

(a3) By (e2) and C3, $U_1\nsubseteq V$
for any $V\in \varphi_{\mathcal{R}_{2}}(2)$.\\
Therefore, $U_1$ satisfies (j1), (j2) and (j3) in Algorithm \ref{a1}.
This follows that $\mathcal{R}_{2}\cup \{U, U_1\}$ is an atomistic lattice which is a 2 order normal subset lattice of $\mathcal{S}$. Analogous to the proof of Case 1,
we can obtain an atomistic lattice $\mathcal{R}_{2}\cup [\emptyset, U]_{\mathcal{Q}}$
which is a 3 order normal subset lattice of $\mathcal{Q}$.

From Cases 1 and 2, we always obtain an atomistic lattice $\mathcal{R}_{2}\cup [\emptyset, U]_{\mathcal{Q}}$
which is a 3 order normal subset lattice of $\mathcal{Q}$.

Continuing as above, we can obtain an atomistic lattice
$$\mathcal{R}_3=\mathcal{R}_2 \cup \bigcup\{[\emptyset, M]_\mathcal{Q}: \ell_{\mathcal{R}_2}(M)=4\}.$$
Obviously, $$\mathcal{R}_3=\mathcal{R}_2 \cup \bigcup\{[\emptyset, M]_\mathcal{Q}: \ell_{\mathcal{Q}}(M)=4, M\in \mathcal{R}_2\}$$ by (e2) in Definition \ref{d41}. Therefore, $\mathcal{R}_3$ is a 4 order normal subset lattice
of $\mathcal{Q}$, and for any $G\in \mathcal{R}_{3}$ with
$\ell_{\mathcal{R}_{3}}(G)\leq 4$ there is no element $H\subsetneq G$ such that $H$ satisfies
(i1), (i2) and (i3) in Algorithm \ref{a1}.

(III). Repeating the preceding proof, we finally obtain an atomistic lattice
$$\mathcal{R}_{m-1}=\mathcal{R}_{m-2} \cup \bigcup\{[\emptyset, M]_\mathcal{Q}: \ell_{\mathcal{Q}}(M)=m, M\in \mathcal{R}_{m-2}\},$$ and for any $W\in \mathcal{R}_{m-1}$ with
$\ell_{\mathcal{R}_{m-1}}(W)\leq m$ there is no element $Z \subsetneq W$ such that $Z$ satisfies
(i1), (i2) and (i3) in Algorithm \ref{a1}. Consequently, by $\ell_{\mathcal{Q}}(M)=m$, we know
that $[\emptyset, M]_\mathcal{Q}=\mathcal{Q}=\mathcal{R}_{m-1}$, and $\mathcal{Q}$ is an output of Algorithm \ref{a1}, completing the proof.
\endproof

Notable that Lemmas \ref{t32}, \ref{t33}, \ref{t34} and \ref{t35} deduce that we can construct all the finite extending cover-preserving geometric lattices of $L$
with the same length by Algorithm \ref{a1}. However, applying the method suggested by G. Cz\'{e}dli and E. T. Schmidt
 in \cite{Czedli10} to the $L$ as depicted in Fig.10, one can only obtain the finite extending cover-preserving geometric lattice as is shown by Fig.15.

 \par\noindent\vskip50pt
\begin{minipage}{11pc}
\setlength{\unitlength}{0.75pt}\begin{picture}(600,410)

\put(300,40){\circle{6}}\put(295,28){\makebox(0,0)[l]
{\footnotesize $\emptyset$}}
\put(300,160){\circle{6}}\put(275,158){\makebox(0,0)[l]
{\footnotesize $\{3\}$}}
\put(300,400){\circle{6}}\put(270,415){\makebox(0,0)[l]
{\footnotesize $\{1,2,3,4,5\}$}}
\put(420,280){\circle{6}}\put(423,280){\makebox(0,0)[l]
{\footnotesize $\{3,5\}$}}
\put(480,280){\circle{6}}\put(483,280){\makebox(0,0)[l]
{\footnotesize $\{2,5\}$}}
\put(540,280){\circle{6}}\put(543,280){\makebox(0,0)[l]
{\footnotesize $\{4,5\}$}}
\put(420,160){\circle{6}}\put(424,150){\makebox(0,0)[l]
{\footnotesize $\{5\}$}}
\put(180,280){\circle{6}}\put(145,280){\makebox(0,0)[l]
{\footnotesize $\{1,3\}$}}
\put(120,280){\circle{6}}\put(85,280){\makebox(0,0)[l]
{\footnotesize $\{1,5\}$}}
\put(60,280){\circle{6}}\put(25,280){\makebox(0,0)[l]
{\footnotesize $\{1,4\}$}}
\put(0,280){\circle{6}}\put(-17,265){\makebox(0,0)[l]
{\footnotesize $\{1,2\}$}}
\put(180,160){\circle{6}}\put(156,150){\makebox(0,0)[l]
{\footnotesize $\{1\}$}}
\put(360,280){\circle{6}}\put(364,280){\makebox(0,0)[l]
{\footnotesize $\{3, 4\}$}}
\put(360,160){\circle{6}}\put(364,160){\makebox(0,0)[l]
{\footnotesize $\{4\}$}}
\put(240,280){\circle{6}}\put(205,280){\makebox(0,0)[l]
{\footnotesize $\{2, 3\}$}}
\put(240,160){\circle{6}}\put(215,160){\makebox(0,0)[l]
{\footnotesize $\{2\}$}}
\put(300,280){\circle{6}}\put(303,280){\makebox(0,0)[l]
{\footnotesize $\{2,4\}$}}

\put(300,43){\line(0,1){114}}
\put(302,162){\line(1,2){57}}
\put(302,42){\line(1,2){57}}
\put(298,42){\line(-1,2){57}}
\put(302,42){\line(1,1){115}}
\put(298,42){\line(-1,1){115}}
\put(298,162){\line(-1,2){57}}
\put(360,163){\line(0,1){114}}
\put(420,163){\line(0,1){114}}
\put(240,163){\line(0,1){114}}
\put(180,163){\line(0,1){114}}
\put(242,162){\line(1,2){57}}
\put(358,162){\line(-1,2){57}}
\put(302,162){\line(1,1){115}}
\put(298,162){\line(-1,1){115}}
\put(422,162){\line(1,2){57}}
\put(422,162){\line(1,1){115}}
\put(178,162){\line(-1,2){57}}
\put(178,162){\line(-1,1){115}}
\put(178,162){\line(-3,2){175}}
\put(358,162){\line(-5,2){295}}
\put(238,162){\line(-2,1){235}}
\put(418,162){\line(-5,2){295}}
\put(242,162){\line(2,1){235}}
\put(362,162){\line(3,2){175}}
\put(300,283){\line(0,1){114}}
\put(242,282){\line(1,2){57}}
\put(182,282){\line(1,1){115}}
\put(122,282){\line(3,2){175}}
\put(62,282){\line(2,1){235}}
\put(2,282){\line(5,2){295}}
\put(358,282){\line(-1,2){57}}
\put(418,282){\line(-1,1){115}}
\put(478,282){\line(-3,2){175}}
\put(538,282){\line(-2,1){235}}
\put(290,5){$G$}
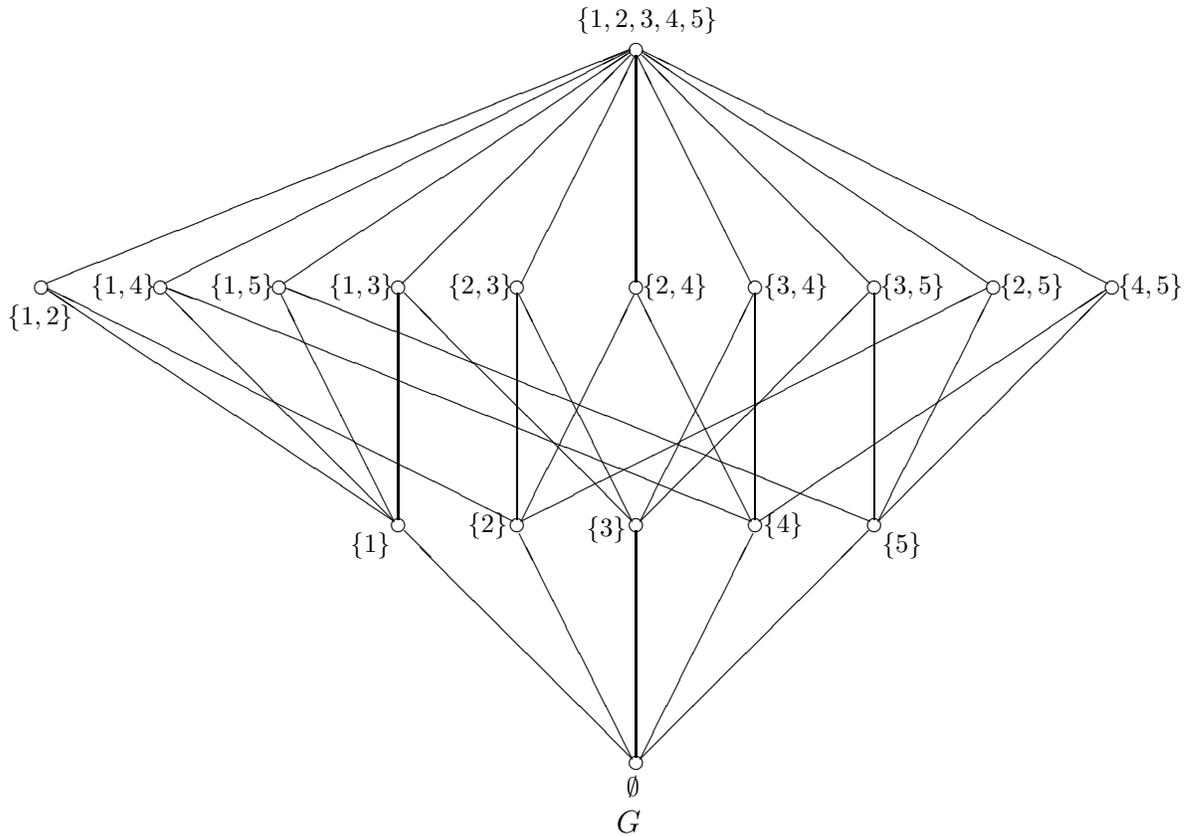
\put(160,-20){ Fig.15 The geometric lattice $G$.}
\end{picture}
\end{minipage}
$$\quad$$

\section{The best geometric lattices}
In this section, we shall construct all the best extending cover-preserving geometric lattices of
$L$. Denote $\overline{\mathfrak{S}}_{k}=\{\mathcal{S}\in \overline{\mathfrak{S}}: |A(\mathcal{S})|=k\}$ for any integer $k>0$. Then we have the following Lemma.

\begin{lemma}\label{t51}
Let $\mathcal{K}\in \overline{\mathfrak{S}}_{k}$ with $k> |J(L)|$. Then there exists an element $\mathcal{H}\in \overline{\mathfrak{S}}_{k-1}$ such that $|\mathcal{H}|< |\mathcal{K}|$.
\end{lemma}
\proof
Since $\mathcal{K}\in \overline{\mathfrak{S}}_{k}$, $L\hookrightarrow^{\prec} \mathcal{K}$ by Lemma \ref{t33}.
Then there exists a lattice
$\mathcal{T}\subseteq \mathcal{K}$ such that
$L\cong \mathcal{T}\hookrightarrow^{\prec} \mathcal{K}$.
Hence there exists a lattice $Q\in \mathfrak{E}(L)$
such that $\mathcal{S}_{Q}=\mathcal{T}\cup A(\mathcal{K})$. This follows from formula (\ref{e31}) that there exists a lattice $\mathcal{T}^{Q}_L$
such that $\mathcal{T}^{Q}_L\subseteq \mathcal{S}_Q$ and
$L\cong \mathcal{T}^{Q}_L \hookrightarrow^{\prec }\mathcal{S}_Q$.
By $k> |J(L)|$, we know that there exists an element $r\in \Delta_Q(L)$
such that $Q-\{r\}\in \mathfrak{E}(L)$. Set $P=Q-\{r\}$ and $R=\{r\}$. Then
$L\cong \mathcal{T}^{P}_{L}\hookrightarrow^{\prec} \mathcal{S}_{P}$ and
\begin{equation}\label{equa111}\mathcal{T}^{P}_{L}=\{F-R: F\in \mathcal{T}^{Q}_L\}\end{equation}
by Lemma \ref{lemma}, and
there exists a set $\sigma \in \mathfrak{I}_{\mathcal{S}_{Q}}(X)$ such that
$\bigcup \sigma\subseteq X-R$ for any $X\in \mathcal{T}^{Q}_L$.
Hence,  as $\mathcal{T}^{Q}_L\subseteq \mathcal{S}_{Q}$, there exists a set $\sigma \in \mathfrak{I}_{\mathcal{K}}(X)$ such
that \begin{equation}\label{000h3}\bigcup \sigma\subseteq X-R\end{equation} for any $X\in \mathcal{T}^{Q}_L$ by
Lemma \ref{l31}. Note that $\mathcal{K}$ is a finite geometric lattice. Then by Lemma \ref{l22}, we have that
\begin{equation}\label{equa112}
\overline{E_{\mathcal{K}}}= \overline{(E-R)_{\mathcal{K}}} =\bigvee_{\mathcal{K}} \sigma
\end{equation} whenever $E\in \mathcal{T}^{Q}_L$, $\sigma \in \mathfrak{I}_{\mathcal{K}}(E)$ and $\bigcup \sigma\subseteq E-R$.

Set
\begin{equation}\label{eq001}
\mathcal{H}=\{W-R: W\in \mathcal{K}, \overline{(W-R)_{\mathcal{K}}}=\overline{W_{\mathcal{K}}}\}.
\end{equation}
Then, by formulas (\ref{equa111}), (\ref{equa112}) and (\ref{eq001}), we know that
\begin{equation}\label{equa113}
\mathcal{T}^{P}_L\subseteq \mathcal{H}.
\end{equation}
Now, we shall show that $L\hookrightarrow^{\prec} \mathcal{H}$ and $\mathcal{H}$ is
a geometric lattice. The proof is made in three steps.

A. $\mathcal{H}$ is a finite atomistic lattice.

From (\ref{eq001}), it is clear that $\mathcal{H}$
is a finite atomistic partially ordered set. Thus, it suffices to show that $\mathcal{H}$ is a lattice.

Suppose $M, N\in \mathcal{H}$. Obviously,
\begin{equation}\label{equationa}R\nsubseteq M \mbox{ and } R\nsubseteq N.\end{equation}
If $M\nparallel N$, then $M\wedge_{\mathcal{H}} N=M\cap N\in \mathcal{H}$.
Now, suppose that $M\parallel N$ and denote $\mathcal{D}_{M, N}=\{G\in \mathcal{H}: G\subseteq M\cap N\}$. There are three cases.

Case i. If $M, N\in \mathcal{K}$, then $M\cap N \in \mathcal{K}$. Clearly, $\overline{(M\cap N-R)_{\mathcal{K}}} =\overline{(M\cap N)_{\mathcal{K}}}=M\cap N$ by (\ref{equationa}).
From (\ref{eq001}), $M\cap N\in \mathcal{H}$. Therefore,
$M\wedge_{\mathcal{H}} N=M\cap N$.

Case j. If $M\in \mathcal{K}$ and $N\notin \mathcal{K}$, then, clearly, $N\cup R\in \mathcal{K}$ since $N\in \mathcal{H}$. Thus $M\cap (N\cup R)\in \mathcal{K}$. By formula (\ref{equationa}), $M\cap (N\cup R)=M\cap N=M\cap N-R$, so that $\overline{[M\cap (N\cup R)]_{\mathcal{K}}} =\overline{(M\cap N)_{\mathcal{K}}}=\overline{(M\cap N-R)_{\mathcal{K}}}$. Hence, $M\cap N\in \mathcal{H}$ by (\ref{eq001}).
Therefore, $M\wedge_{\mathcal{H}} N=M\cap N$.

Case k. If $M, N\notin \mathcal{K}$, then similar to the proof of Case j,
we have $M\cup R,N\cup R\in \mathcal{K}$. Then $(M\cap N)\cup R\in \mathcal{K}$.
There are two subcases.

Subcase $1^\circ$. If $\overline{[(M\cap N)\cup R]_{\mathcal{K}}} =\overline{(M\cap N)_{\mathcal{K}}}$,
then by (\ref{eq001}), $M\cap N\in \mathcal{H}$ which is the maximum
element of $\mathcal{D}_{M, N}$. Thus, $M\wedge_{\mathcal{H}} N=M\cap N$.

Subcase $2^\circ$. If$\overline{[(M\cap N)\cup R]_{\mathcal{K}}}\neq \overline{(M\cap N)_{\mathcal{K}}}$,
then $\overline{(M\cap N)_{\mathcal{K}}}\subsetneq \overline{[(M\cap N)\cup R]_{\mathcal{K}}}$, and $\overline{[(M\cap N)\cup R]_{\mathcal{K}}}= (M\cap N)\cup R$ since $(M\cap N)\cup R\in \mathcal{K}$.
Thus $$M\cap N \subseteq \overline{(M\cap N)_{\mathcal{K}}}
\subsetneq \overline{[(M\cap N)\cup R]_{\mathcal{K}}}= (M\cap N)\cup R,$$ which means that
$M\cap N =\overline{(M\cap N)_{\mathcal{K}}}\in \mathcal{K}$. Then, by (\ref{equationa}), we know that
$M\cap N=\overline{(M\cap N)_{\mathcal{K}}}=\overline{[(M\cap N)-R]_{\mathcal{K}}}$. Therefore,
$M\cap N\in \mathcal{H}$ by (\ref{eq001}), it follows that $M\wedge_{\mathcal{H}} N=M\cap N$.

Subcases $1^\circ$ and $2^\circ$ imply that $M\wedge_{\mathcal{H}} N=M\cap N$ if $M, N\notin \mathcal{K}$.

Therefore, from Cases i, j and k, $\mathcal{H}$ is a finite atomistic lattice.

B. $\mathcal{H}$ is a finite geometric lattice.

Inasmuch as we possess A it suffices to show that
$\mathcal{H}$ is a semimodular lattice. Let $M, N\in \mathcal{H}$ and $M, N\succ_{\mathcal{H}} M\cap N$.
Obviously, $M\parallel N$. Next, we shall prove that
$M\vee_{\mathcal{H}}N\succ_{\mathcal{H}} M, N$. There are three cases as follows.

Case a. If $M, N\in \mathcal{K}$, then $M\cap N\in \mathcal{K}$.
By formulas (\ref{eq001}) and (\ref{equationa}), $[\emptyset, M]_{\mathcal{H}}=[\emptyset, M]_{\mathcal{K}}$ and
$[\emptyset, N]_{\mathcal{H}}=[\emptyset, N]_{\mathcal{K}}$.
Thus $M, N\succ_{\mathcal{K}} M\cap N$, which together with $\mathcal{K}$ is a geometric lattice yields that  $M\vee_{\mathcal{K}} N\succ_{\mathcal{K}} M, N$.

We claim that $$\overline{(M\vee_{\mathcal{K}} N)_{\mathcal{K}}}=\overline{(M\vee_{\mathcal{K}} N-R)_{\mathcal{K}}}.$$
Otherwise, $\overline{(M\vee_{\mathcal{K}} N-R)_{\mathcal{K}}}\subsetneq \overline{(M\vee_{\mathcal{K}} N)_{\mathcal{K}}}$. It is clear that
$M\subsetneq M\cup N\subseteq \overline{(M\vee_{\mathcal{K}} N-R)_{\mathcal{K}}}$ by (\ref{equationa}).
Note that $\overline{(M\vee_{\mathcal{K}} N)_{\mathcal{K}}}=M\vee_{\mathcal{K}} N$ since
$M\vee_{\mathcal{K}} N\in \mathcal{K} $.
Thus $$M\subsetneq \overline{(M\vee_{\mathcal{K}} N-R)_{\mathcal{K}}}\subsetneq \overline{(M\vee_{\mathcal{K}} N)_{\mathcal{K}}}=M\vee_{\mathcal{K}} N,$$ contrary to the fact that
$M\vee_{\mathcal{K}} N\succ_{\mathcal{K}} M$. Therefore,  $\overline{(M\vee_{\mathcal{K}} N)_{\mathcal{K}}}=\overline{(M\vee_{\mathcal{K}} N-R)_{\mathcal{K}}}$, and then by (\ref{eq001}),
$M\vee_{\mathcal{K}} N-R\in \mathcal{H}$. Hence, the condition $M\vee_{\mathcal{K}} N\succ_{\mathcal{K}} M, N$
deduces that $$M\vee_{\mathcal{K}} N-R=M\vee_{\mathcal{H}} N\succ_{\mathcal{H}} M, N.$$

Case b. If $M, N \notin \mathcal{K}$, then $M\cup R, N\cup R, (M\cap N)\cup R\in \mathcal{K}$ since $M, N, M\cap N \in \mathcal{H}$. Thus
\begin{equation}\label{0001}
M\cup R=\overline{(M\cup R)_{\mathcal{K}}}=\overline{M_{\mathcal{K}}},  N\cup R=\overline{(N\cup R)_{\mathcal{K}}}=\overline{N_{\mathcal{K}}}
\end{equation} and \begin{equation}\label{000h1}(M\cap N)\cup R=\overline{[(M\cap N)\cup R]_{\mathcal{K}}}=\overline{(M\cap N)_{\mathcal{K}}}\end{equation} by (\ref{eq001}).

On the other hand, we claim that $M\cup R\succ_{\mathcal{K}} (M\cap N)\cup R$. Otherwise,
there exists an atom $I\subseteq M-N$ such that $M\cup R\supsetneq \overline{[(M\cap N)\cup R\cup I]_{\mathcal{K}}}\supsetneq(M\cap N)\cup R$.
Then by formula (\ref{000h1}),
$\overline{[(M\cap N)\cup R\cup I]_{\mathcal{K}}}=\overline{[(M\cap N)\cup I]_{\mathcal{K}}}$,
which means that $\overline{[(M\cap N)\cup I]_{\mathcal{K}}}-R\in \mathcal{H}$.
Thus $M\supsetneq \overline{[(M\cap N)\cup I]_{\mathcal{K}}}-R\supsetneq M\cap N$,
contrary to the fact that $M\succ_{\mathcal{H}} M\cap N$.
Therefore, \begin{equation}\label{000h2}M\cup R\succ_{\mathcal{K}} (M\cap N)\cup R.\end{equation}
Similarly, we have $N\cup R\succ_{\mathcal{K}} (M\cap N)\cup R$. Thus $(M\cup R) \vee_{\mathcal{K}} (N\cup R)\succ_{\mathcal{K}} M\cup R, N\cup R$, and which means that
\begin{equation}\label{0002}
\overline{(M\cup N\cup R)_{\mathcal{K}}}= (M\cup R) \vee_{\mathcal{K}}
(N\cup R) \succ_{\mathcal{K}} M\cup R, N\cup R.
\end{equation}
 We claim that
$\overline{[\overline{(M\cup N\cup R)_{\mathcal{K}}}]_{\mathcal{K}}}=
\overline{[\overline{(M\cup N\cup R)_{\mathcal{K}}}-R]_{\mathcal{K}}}$.
Otherwise, $$\overline{(M\cup N\cup R)_{\mathcal{K}}}\supsetneq \overline{[\overline{(M\cup N\cup R)_{\mathcal{K}}}-R]_{\mathcal{K}}}\supseteq \overline{(M \cup N)_{\mathcal{K}}}$$ since
$\overline{[\overline{(M\cup N\cup R)_{\mathcal{K}}}]_{\mathcal{K}}}=\overline{(M\cup N\cup R)_{\mathcal{K}}}$
and $\overline{(M \cup N\cup R)_{\mathcal{K}}} -R\supseteq M\cup N$.
However, $\overline{(M \cup N)_{\mathcal{K}}}\supseteq \overline{M_{\mathcal{K}}}=M\cup R$ by
(\ref{0001}). Hence $\overline{(M \cup N)_{\mathcal{K}}}=\overline{(M \cup N\cup R)_{\mathcal{K}}}$,
a contradiction. Therefore, $\overline{[\overline{(M\cup N\cup R)_{\mathcal{K}}}]_{\mathcal{K}}}=
\overline{[\overline{(M\cup N\cup R)_{\mathcal{K}}}-R]_{\mathcal{K}}}$. This follows that $\overline{(M\cup N\cup R)_{\mathcal{K}}}-R\in \mathcal{H}$ by (\ref{eq001}).
Further, by formula (\ref{0002}), $$\overline{(M\cup N\cup R)_{\mathcal{K}}}-R=\overline{(M\cup N)_{\mathcal{K}}}-R=M\vee_{\mathcal{H}} N\succ_{\mathcal{H}}M, N.$$

Case c. If $M\notin \mathcal{K}$ and $N \in \mathcal{K}$, or $M\in \mathcal{K}$ and $N \notin \mathcal{K}$, say, $M\notin \mathcal{K}$ and $N \in \mathcal{K}$, then $M\cup R, N, M\cap N\in \mathcal{K}$ since $M , N, M\cap N\in \mathcal{H}$.
Hence,
\begin{equation}\label{0003}
M\cup R=\overline{(M\cup R)_{\mathcal{K}}}=\overline{M_{\mathcal{K}}}
\end{equation} by (\ref{eq001}).

Similar to the proof of formula (\ref{000h2}), we have that $(M\cup R)
\succ_{\mathcal{K}} M\cap N$. On the other hand, similar to the proof of Case a, we know that $N\succ_{\mathcal{K}} M\cap N$. Thus $(M\cup R) \vee_{\mathcal{K}} N\succ_{\mathcal{K}} M\cup R, N$,
and which means that
\begin{equation}\label{0004}
\overline{(M\cup N\cup R)_{\mathcal{K}}}= (M\cup R) \vee_{\mathcal{K}}
N \succ_{\mathcal{K}} M\cup R, N.
\end{equation}

Analogous to the proof of Case b, we know that
$$\overline{(M\cup N\cup R)_{\mathcal{K}}}-R=\overline{(M\cup N)_{\mathcal{K}}}-R=M\vee_{\mathcal{H}} N\succ_{\mathcal{H}} M, N$$
by formulas (\ref{eq001}), (\ref{0003}) and (\ref{0004}).

In summary, $\mathcal{H}$ is a finite geometric lattice.

C. $\mathcal{T}^{P}_L\hookrightarrow^\prec \mathcal{H}$.

Let $M, N\in \mathcal{T}^{P}_L$.
Then there are two elements
$E, F\in \mathcal{T}^{Q}_L\subseteq \mathcal{S}_Q$ such that $M=F-R$ and $N=E-R$ by (\ref{equa111}).
Thus by $L\cong \mathcal{T}^{P}_L\cong \mathcal{T}^{Q}_L$ and (\ref{equa111}),
\begin{equation}\label{ww1}
 M\vee_{\mathcal{T}^{P}_L} N =(F\vee_{\mathcal{T}^{Q}_L} E)-R.
\end{equation}
As $F, E\in \mathcal{T}^{Q}_L$, we have that $$F\vee_{\mathcal{K}} E=\overline{(M\cup N)_{\mathcal{K}}}=\overline{[\overline{(M\cup N)_{\mathcal{K}}}-R]_{\mathcal{K}}}$$ by formula (\ref{000h3}) and (\ref{equa112}).
Hence, by formulas (\ref{eq001}), \begin{equation}\label{ww2}M\vee_{\mathcal{H}} N=\overline{(M\cup N)_{\mathcal{K}}}-R=F\vee_{\mathcal{K}} E-R.\end{equation}

Clearly, by formula (\ref{wwwb}),
we know that $$F\vee_{\mathcal{T}^{Q}_L} E=F\vee_{\mathcal{K}} E,$$
which together with formulas (\ref{ww1}) and (\ref{ww2}) clearly leads to $M\vee_{\mathcal{H}} N =M\vee_{\mathcal{T}^{P}_L} N$. Therefore, $\mathcal{T}^{P}_L \hookrightarrow^{\vee} \mathcal{H}$.
On the other hand, by formula (\ref{equa111}), we know that
$M\wedge_{\mathcal{T}^{P}_L} N=M\cap N$ for
any $M, N\in \mathcal{T}^{P}_L$.
Thus $\mathcal{T}^{P}_L \hookrightarrow^{\wedge} \mathcal{H}$ by the proof of A.
Hence, $\mathcal{T}^{P}_L \hookrightarrow \mathcal{H}$.

Obviously, $L\cong \mathcal{T}^{P}_L\cong \mathcal{T}^{Q}_L$ together with
(\ref{eq001}) means
that $\ell(\mathcal{T}^{P}_L)=\ell(\mathcal{H})$.
Therefore, $\mathcal{T}^{P}_L \hookrightarrow^{\prec} \mathcal{H}$ since $\mathcal{T}^{P}_L \hookrightarrow \mathcal{H}$, completing the proof of C.

Finally, from B, C and Lemma \ref{t35}, we know
that $\mathcal{H}\in \overline{\mathfrak{S}}_{k-1}$
and $|\mathcal{H}|< |\mathcal{K}|$. This completes the proof.
\endproof

Let $G$ be a finite geometric lattice. It is clear that if $L\hookrightarrow^{\prec} G$, then
there exists a sublattice $[x, y]$ of $G$ with $\ell([x, y])=\ell(L)$ such that $L\hookrightarrow^{\prec} [x, y]$. Clearly, $[x, y]$ is also a geometric lattice. Therefore, by Lemmas \ref{t32}, \ref{t33}, \ref{t34}, \ref{t35} and \ref{t51}, we have the following theorem.
\begin{theorem}\label{t52}
Every best extending cover-preserving geometric lattice of $L$
is the best one in $\overline{\mathfrak{S}}_{|J(L)|}$.
\end{theorem}

\begin{example}\label{ex51}
\emph{Consider the lattice $L$ in Example \ref{ex31} again.
If $U=\{1,2,4,5\}$ in Step 1 of Algorithm \ref{a1}.
Then $\mathcal{Q}=\mathcal{S}_P\cup \{U\}$ is the output lattice of Algorithm \ref{a1}
(the lattice $\mathcal{Q}$ as represented in Fig.16).}
\end{example}
\par\noindent\vskip50pt
\begin{minipage}{11pc}
\setlength{\unitlength}{0.75pt}\begin{picture}(600,250)

\put(300,40){\circle{6}}\put(295,28){\makebox(0,0)[l]
{\footnotesize $\emptyset$}}
\put(300,120){\circle{6}}\put(275,118){\makebox(0,0)[l]
{\footnotesize $\{3\}$}}
\put(300,280){\circle{6}}\put(270,295){\makebox(0,0)[l]
{\footnotesize $\{1,2,3,4,5\}$}}
\put(460,200){\circle{6}}\put(463,198){\makebox(0,0)[l]
{\footnotesize $\{3,5\}$}}
\put(460,120){\circle{6}}\put(463,118){\makebox(0,0)[l]
{\footnotesize $\{5\}$}}
\put(140,200){\circle{6}}\put(104,198){\makebox(0,0)[l]
{\footnotesize $\{1,3\}$}}
\put(140,120){\circle{6}}\put(114,118){\makebox(0,0)[l]
{\footnotesize $\{1\}$}}
\put(380,200){\circle{6}}\put(385,200){\makebox(0,0)[l]
{\footnotesize $\{3, 4\}$}}
\put(380,120){\circle{6}}\put(385,120){\makebox(0,0)[l]
{\footnotesize $\{4\}$}}
\put(220,200){\circle{6}}\put(185,198){\makebox(0,0)[l]
{\footnotesize $\{2, 3\}$}}
\put(220,120){\circle{6}}\put(196,120){\makebox(0,0)[l]
{\footnotesize $\{2\}$}}

\put(300,200){\circle{6}}\put(304,205){\makebox(0,0)[l]
{\footnotesize $\{1,2,4,5\}$}}

\put(380,123){\line(0,1){74}}
\put(220,123){\line(0,1){74}}
\put(140,123){\line(0,1){74}}
\put(460,123){\line(0,1){74}}

\put(300,43){\line(0,1){74}}
\put(302,122){\line(1,1){76}}
\put(298,122){\line(-1,1){76}}

\put(222,202){\line(1,1){76}}
\put(378,202){\line(-1,1){76}}
\put(302,122){\line(2,1){155}}
\put(298,122){\line(-2,1){155}}
\put(142,202){\line(2,1){155}}
\put(458,202){\line(-2,1){155}}

\put(302,42){\line(2,1){155}}
\put(298,42){\line(-2,1){155}}
\put(302,42){\line(1,1){76}}
\put(298,42){\line(-1,1){76}}
\put(378,122){\line(-1,1){76}}
\put(222,122){\line(1,1){76}}
\put(142,122){\line(2,1){155}}
\put(300,203){\line(0,1){74}}
\put(458,122){\line(-2,1){155}}
\put(290,0){$\mathcal{Q}$}
\put(230,-25){Fig.16 The lattice $\mathcal{Q}$.}

\end{picture}
\end{minipage}
$$\quad$$
\\
Obviously, $\mathcal{Q}\in \overline{\mathfrak{S}}_{5}$. Further,
we know that $\mathcal{Q}$ is the unique
best extending cover-preserving geometric lattice of $L$ in the sense of isomorphism.

\section{Conclusions}
In this paper, we proposed an algorithm to calculate all the
best extending cover-preserving geometric lattice $G$ of a given semimodular lattice $L$ and proved that $|A(G)|=|J(L)|$ and $\ell(G)=\ell(L)$.
It is worth pointing out that every different $U$ (resp. $W$) in Algorithm \ref{a1} leads to a different output, and the computational complexity of Algorithm \ref{a1} is likely to grow rapidly as $|J(L)|$ and $\ell(L)$ grow.

\section*{Data availability statements}
The datasets generated during and/or analysed during the current study are available from the corresponding author on reasonable request.

\end{document}